\documentclass[12pt,a4paper]{article}
\usepackage[utf8]{inputenc}
\usepackage[english]{babel}
\usepackage{amsmath}
\usepackage{amsfonts}
\usepackage{amssymb}
\usepackage{makeidx}
\usepackage{graphicx}
\usepackage{lmodern}
\usepackage[left=3cm,right=3cm,top=2.5cm,bottom=4cm]{geometry}
\usepackage[utf8]{inputenc}
\usepackage[T1]{fontenc}
\usepackage{lmodern}
\usepackage{graphicx}
\usepackage[figurename=Fig.,labelfont=bf,labelsep=period]{caption}
\usepackage{subcaption}
\usepackage{amsmath}
\usepackage{adjustbox}
\usepackage{natbib}
\usepackage{titlesec}
\titleformat*{\section}{\normalsize \bfseries}
\titleformat*{\subsection}{\normalsize \normalshape}
\titleformat*{\subsubsection}{\normalsize \itshape}
\titleformat*{\paragraph}{\normalsize \itshape}
\titleformat*{\subparagraph}{\normalsize\itshape}
\usepackage[colorlinks=true,citecolor=blue,linkcolor=black]{hyperref}
\usepackage{graphicx} 
\usepackage{caption}
\graphicspath{{Figures/}}
\title{Scaling laws for the rigid-body response of masonry structures under blast loads}
\author{Filippo Masi$^a$ \and Ioannis Stefanou$^a$ \and Victor Maffi-Berthier$^b$}
\date{\small
    $^a$Institut de Recherche en Génie Civil et Mécanique,\\
    UMR 6183, CNRS, École Centrale de Nantes, Université de Nantes, 1 rue de la N\"{o}e, F-44300, Nantes, France.\\
    Email: \href{mailto:filippo.masi@ec-nantes.fr}{filippo.masi@ec-nantes.fr}, \href{mailto:ioannis.stefanou@ec-nantes.fr}{ioannis.stefanou@ec-nantes.fr}\\%
    $^b$Ing\'{e}rop Conseil et Ing\'{e}nierie,\\
    18 rue des Deux Gares, F-92500, Rueil-Malmaison, France.\\
    Email: \href{mailto:victor.maffi-berthier@ingerop.com}{victor.maffi-berthier@ingerop.com}\\[2ex]%
    \today
}

\begin{document}
\maketitle
\normalshape
\begin{abstract}
\noindent The response of masonry structures to explosions can be hardly investigated relying only on numerical and analytical tools. 
Experimental tests are of paramount importance for improving the current comprehension and validate existing models. However, experiments involving blast scenarios are, at present, partial and limited in number, compared to tests under different dynamic conditions, such as earthquakes. The reason lies on the fact that full-scale blast experiments present many difficulties, mainly due to the nature of the loading action. Experiments in reduced-scale offer instead greater flexibility. Nevertheless, appropriate scaling laws for the response of masonry structures under blast excitations are needed before performing such tests.\\
We propose here new scaling laws for the dynamic, rigid-body response and failure modes of masonry structures under blast loads. This work takes its roots from previous studies, where closed-form solutions for the rocking response of slender blocks due to explosions have been derived and validated against numerical and experimental tests.\\
The proposed scaling laws are here validated with detailed numerical simulations accounting for combined rocking, up-lifting, and sliding mechanisms of monolithic structures. Then, the application to multi-drum stone columns is considered. In particular, we show that, whilst the presence of complex behaviors, such as wobbling and impacts, similarity is assured. The developments demonstrate their applicability in the design of reduced-scale experiments of masonry structures.
\end{abstract}


\section{Introduction}
\label{sec:intro_sim}
\noindent The investigation of the response and strength of masonry structures subjected to blast loading is a major issue in protecting such assets against explosions. In particular, we refer here to stone and brick structures of ancient and classical monuments, which have been repeatedly, and can be in the future, exposed to accidental or deliberate blasts (e.g. the Parthenon in Athens in 1687, the Monumental Arch of Palmyra in 2015, the recent blast in Beirut, on 4 August 2020, or the partial destruction of Ghazanchetsots Cathedral, on 10 October 2020). Consequently, appropriate evaluations and assessment of the consequences of a blast event on masonry buildings and their key load-carrying members are mandatory.

In the last decades, increasing scientific interest has being shown in the blast response and potential damage of masonry structures. Attention has been mostly focused on regular geometries, i.e., on flat walls. This was achieved in terms of both experimental campaigns and numerical modeling, see e.g. \cite{gabrielsen1975response,KEYS2017229,LI2017107,gebbeken,hstamwall}.\\
Research activity for the dynamic response of ancient and classical stone monuments in a blast event is, at the moment, modest and insufficient to draw appropriate preservation strategies and in-situ protections \cite{hstamwall,masi2020discrete}. This is due to the fact that blast testing of masonry structures is particularly challenging due to its mechanical properties.\\
Furthermore, such kind of assets are usually characterized by a \textit{non-standard} geometry. The term \textit{non-standard} geometry denotes here those assets being composed of arches, vaults, domes, and any other structural element that is substantially different from planar walls \citep{cnrs,forgacs2017minimum,forgacs2018influence,VANNUCCI20192,stockdale2020seismic,kassotakis2020three,masithesis}. The complexity stemming from non-standard geometries adds non-trivial difficulties into the design of experimental tests.\\

Some important in-situ tests of masonry exist \citep[][, among others]{gabrielsen1975response,varma,abou,KEYS2017229,gilbert2002performance,LI2017107,gebbeken}. These tests provided some information but reproducting the same conditions is highly demanding due to associated high cost, safety issues, reduced repeatability, technical complications, etc. \cite{draganic2018overview}. For instance, explosive charge weights need to be as small as possible for safety reasons. Moreover, specialized personnel is required to built masonry and the execution of repeated, in-situ experiments is hardly possible. It is needless to mention that performing real-scale tests in proper laboratory conditions is impossible unless large investments are made.\\

Conversely, reduced-scale tests \citep[e.g.][]{wang2012experimental} offer many advantages, such as reduced cost, reduced hazard and risk associated to the safety of the testing environment and of the personnel. Nevertheless, the design of small-scale tests requires appropriate scaling laws in order to guarantee similarity. Similitude theory provides the conditions to design a scaled model of a prototype (full-scale structure) and to predict the structural response of the prototype from the scaled results \cite[cf. ][]{CAKTI2016224}. Similarity of the blast wave propagation and load need to be assured between the scaled model/system and the prototype. The same holds true for the specimen where scaling laws for the dynamic structural response have to be derived to assure similarity between the model and the prototype. Notice, though, that deriving adequate scaling laws for both the blast load and the structure is not trivial due to the different underlying physics.\\
As far it concerns blast loads, the Hopkinson-Cranz scaling law is widely used \cite{hopki,cranz}. This scaling law assumes similarity at constant scaled distances, $Z=R/W^{1/3}$, with $W$ the TNT equivalent explosive weight and $R$ the stand-off distance between the charge and the target. Accordingly, the scaled distance in the model needs to be equal to the scaled distance in the prototype. As a result, the overpressure peak, originated by the explosion, remains the same in the prototype and the small-scale system. The Hopkinson-Cranz scaling law has been successfully applied in numerous applications. For more details, we refer to \citet{199173}. Nevertheless, it is worth emphasizing that a scaling preserving the overpressure peak is disadvantageous and inadequate for performing reduced-scale tests in a laboratory environment. Indeed, most of the benefits of reduced-scale testing with respect to full-scale testing are lost, as the intensity of the explosive load is un-scaled. As a result, the scaling laws available in the current literature are not a viable mean to design safe, laboratory tests of whole masonry structures.\\

The main goal of this paper is to provide and test a new set of scaling laws for masonry structures subjected to explosions and guide towards experimental testing of non-standard assets. In particular, based on the developments of the rigid-body response of monolithic structures against blast loads \cite{rockingmasi}, we derive here scaling laws for masonry structures in a blast event. This is done considering the rigid-body response, occurring after material failure. Indeed, under the action of the impulsive loading arising from an explosion, local failure of the joints and blocks results, in most cases, in failure modes which can be described, depending on the kinematics of the failure mode, as a rigid-body motion \citep[cf.][]{masithesis}. Nevertheless, both the rigid-body and the material responses can be scaled according to the approach here proposed \citep[cf.][]{masithesis}.\\

This paper is structured as follows. First, we recall the characteristic features of blast waves, the corresponding loading and their modeling. Second, the problem of the rigid-body motion of masonry structures is stated. Then, scaling laws for the rigid-body response are derived and validated for cases of monolithic (rigid) prototypes and models. Finally, we show that the derived scaling laws hold true for deformable blocky-structures, by investigating the response of multi-drum columns: typical examples of key load-carrying elements in ancient masonry structures and monuments.\\
It should be mentioned that this work is a first step towards the design of reduced-scale experiments of masonry structures, providing for the the first time appropriate scaling laws assuring the similarity of both blast loading and structural dynamic response. Further investigations accounting for more detailed characterization of blast wave propagation and secondary effects, such as clearing and eventual blast waves focalization phenomena, as well as different masonry assets are needed.

\section{Blast actions and model}
\label{sec:blast}
\noindent Explosions produce blast waves of high-pressure, accompanying high-temperature and supersonic expansion of gases.
The abrupt increase of the pressure carried by a blast wave can produce severe structural damage. When the primary shock strikes a target, the so-called reflected overpressure, $P_r$, originates. Figure \ref{fig:blast1} shows schematically the time evolution of $P_r$, which is determined by the arrival time of the shock wave, $t_A$, the overpressure peak, $P_{ro}$, the positive phase duration, $t_o$, the negative phase duration, $t_{o-}$, and the underpressure peak, $P_{ro-}$. These parameters are functions of the distance, $R$, for the explosive source and of the explosive weight, $W$ (conventionally expressed in TNT equivalent). Herein we consider only the positive phase of the blast wave \cite[safety approach, see ][]{rockingmasi,masithesis}. 

The simulation of a blast can be conducted by using different approaches \cite{masithesis}. Empirical methods, fitted on actual experiments, are here considered. Such methods rely on best-fit interpolations of experimental results and mainly on those of Kingery and Bulmash \cite{kingery}, which allow to determine the blast parameters and pressure loading from the knowledge of the trinitrotoluene (TNT) equivalent explosive weight, $W$, and the Hopkinson-Cranz scaled distance, $Z=R/\sqrt[3]{W}$ (see Appendix I).\\
The time evolution of the positive phase of the reflected pressure is modeled with the well established \textit{modified Friedlander equation} \cite{friedlander},
\begin{equation}
P_{r}(t)=  P_{ro} \left(1-\dfrac{t}{t_o} \right)\Big(1 - \mathcal{H}[t - t_o] \Big)  \exp{\left(-d \dfrac{t}{t_o} \right)},
\label{eq:fried}
\end{equation}
where $\mathcal{H}[\cdot]$ denotes the Heaviside (step) function, $d$ is the exponential decay coefficient, and $t_A$ is taken as the origin of the time axis.
The impulse $i_{ro}$ associated to the positive phase, which represents the area beneath the pressure curve, reads
\begin{equation}
i_{ro}=\int_{0} ^{t_o} P_{r} \, dt=\left[ e^{-d}+d-1\right] \dfrac{P_{ro}\,t_o }{d^2}.
\label{eq:imp_p}
\end{equation}
The above equation allows to determine the exponential decay coefficient, $d$, by equating it with the best-fit interpolation of $i_{ro}$ from experiments (see Appendix I).
%
%
\subsection{First-order approximation}
\begin{figure*}[ht]
\centering
\begin{subfigure}[t]{0.38\textwidth}
  \centering
  \includegraphics[width=\linewidth]{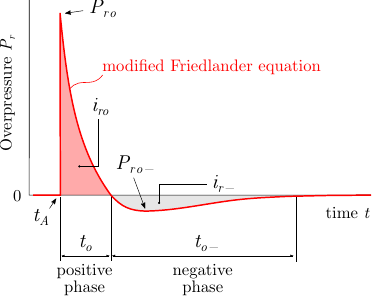}
	\caption{\footnotesize modified Friedlander equation}
	\label{fig:blast1}
\end{subfigure} \hspace{1cm}
\begin{subfigure}[t]{0.38\textwidth}
  \centering
  \includegraphics[width=\linewidth]{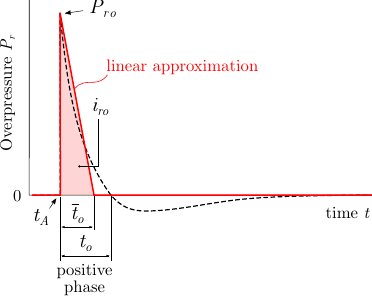}
	\caption{\footnotesize first-order approximation}
	\label{fig:blast2}
\end{subfigure}
\caption{Time evolution of overpressure due to an explosion acting on a target: (a) modified Friedlander equation and (b) first-order approximation.}
		\label{fig:blast}
\end{figure*}

\noindent One way of modeling the blast pressure is with a piece-wise linear function, as is the approach of \citet{KRAUTHAMMER20001}, also recommended in \citet{ufc08}. In this case, the positive phase is approximated as a linearly decaying triangular pulse (see Figure \ref{fig:blast2})
\begin{equation}
P_{r}(t)=  P_{ro} \left(1-\dfrac{t}{\bar{t}_o} \right)\Big(1 - \mathcal{H}[t - t_o] \Big),
\label{eq:friedlin}
\end{equation}
with $\bar{t}_o$ being the linear load duration, set such that the impulse given by the triangular load matches that given by the best fitting interpolations (Appendix I), and namely
\begin{equation}
\bar{t}_o = 2 \dfrac{i_{ro}}{P_{ro}}.
\end{equation}

\section{Problem statement}
\label{sec:similarity}
\begin{figure*}[ht]
\centering
\begin{subfigure}[b]{0.35\textwidth}
  \centering
  \includegraphics[width=\linewidth]{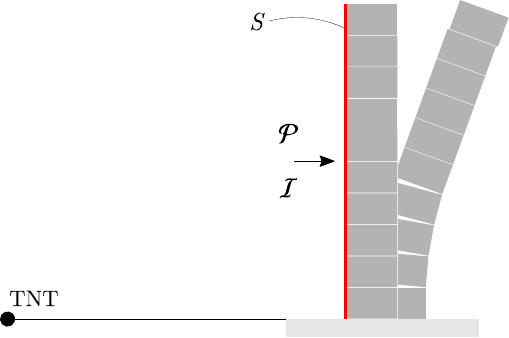}
  \caption{\footnotesize}
\end{subfigure} \hspace{1cm}
\begin{subfigure}[b]{0.35\textwidth}
  \centering
  \includegraphics[width=\linewidth]{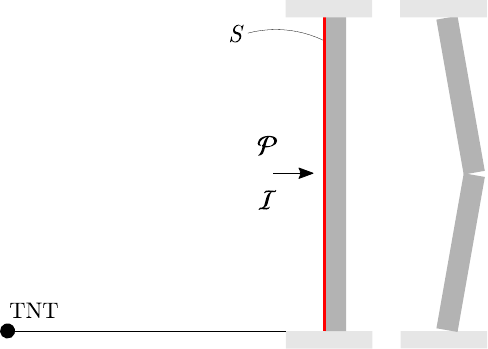}
  \caption{\footnotesize}
\end{subfigure}
\caption{Representative scheme of the problem, e.g. (a) a multidrum column or (b) a one-way spanning wall, under blast loading.}
		\label{fig:scheme}
\end{figure*}

\noindent Let us consider a masonry structure, of arbitrary shape, composed of masonry units, interacting one with the other through interfaces with friction angle $\varphi$, (non-)associative sliding behavior, and zero cohesion and tensile strength \cite{Godio2018}. The structure is subjected to the load of an explosion and undergoes a rigid-body motion, see Figure \ref{fig:scheme}. We assume that deformations are negligible with respect to the rigid-body response, which are developed during material failure. Such a hypothesis is usually realistic for masonry structures under low confinement, e.g. for columns and retaining walls or during rocking-like failure of internal segments. The consequence of such an assumption is discussed further in the paper, by accounting for the elastic deformations of masonry elements.\\
 The loading force is characterized by its maximum specific thrust $\mathcal{P}$ and the maximum specific impulse $\mathcal{I}$. For targets small enough to assume that the blast wave acts simultaneously and uniformly on the impinged surfaces, the maximum specific thrust and impulse can be computed as the overpressure peak and impulse of an equivalent blast load, acting on the center of the mass of the structure. Following the developments in \citet{rockingmasi}, we consider only the pressure load applied on the front surface $S$ (incident surface, Figure \ref{fig:scheme}) of the target. We further assume that the blast wave impinges all points of $S$ at the same time (simultaneously) and with the same magnitude (uniformly). The limitations of these assumptions are further discussed in the applications to multi-drum masonry columns and in Appendix II. Moreover, as stated above, we account only for the positive phase of the blast load. The negative phase is expected to have a stabilizing effect. Therefore, it can be neglected to have safe estimates of the system response, as shown in \citet{rockingmasi}.\\

We denote the geometry of the structure by (i) a characteristic length $l$ (e.g., the height of the structure in Figure \ref{fig:scheme}), (ii) dimensionless length ratios $l_i$, which relate all other lengths to the characteristic one, and (iii) generalized angles $\alpha_i$. The structure is further characterized by mass $m$, mass moment of inertia $J$ about some specific axis, and non-dimensional mass moments of inertia ratios $J_i$, which relate all other components of the rotational inertia tensor to the characteristic mass moment of inertia. We denote the gravitational acceleration with $g$ and the friction coefficient of the interfaces with $\mu= \tan\varphi$. Coulomb friction is adopted. The material has density $\rho$. For each block constituting the structure, we identify the sliding distance, $x$, the linear velocity and acceleration, respectively $\dot{x}$ and $\ddot{x}$, the rocking angle, $\theta$, and the angular velocity and acceleration of the blocks, $\dot{\theta}$ and $\ddot{\theta}$, respectively.\\

From the $\pi$ Theorem \citep{bertrand1878homogeneite}$-$see also \citep{199173}$-$, we identify the following terms for the rigid-body response of the blocks:
\begin{subequations}
	\begin{align}
	\pi_{01}=J_i, \quad \pi_{02}= \dfrac{x}{l}, \quad	\pi_{03}= \theta,\\
	\pi_{11}=\ddot{\theta} \dfrac{l}{g},\quad
	\pi_{12} =\dot{\theta} \sqrt{\dfrac{l}{g}},\quad 
	\pi_{13} = \dfrac{\ddot{x}}{g},\quad
	\pi_{14} = \dfrac{\dot{x}}{\sqrt{lg}},\quad 
	\pi_{15}= t  \sqrt{\dfrac{g}{L}},\\
	\pi_{21}=\mu,\quad
	\pi_{22} =\dfrac{J}{ml^2},\quad
	\pi_{23}= \dfrac{\mathcal{P} l^2 }{m g},\quad
	\pi_{24}= \dfrac{\mathcal{I}}{m} \sqrt{\dfrac{l^3}{g}}	
\end{align}
\end{subequations}
Terms $\pi_{01}-\pi_{03}$ represent the geometric similarity and terms $\pi_{11}-\pi_{15}$ identify the kinematic similarity, i.e., the response of the system in terms of linear and angular displacements, velocities, and accelerations. The remaining four terms $\pi_{21}-\pi_{24}$ determine dynamic similarity in terms of rigid-body motion.\\
By imposing the equivalence of the aforementioned $\pi$ terms, one can investigate and predict the response of a full-scale system (namely, a prototype) by studying the response of a reduced scale system (i.e., a model), satisfying the similarity statements.\\

Let us suppose that both prototype and model consist of blocks with a uniformly distributed mass, $m$ and $\tilde{m}$, respectively. The geometric scaling performed on the model is determined by the ratio of the characteristic lengths of the two systems, i.e.,
\begin{equation*}
\lambda= \dfrac{\tilde{l}}{l},
\end{equation*}
with superscript `$\sim$' denoting model's quantities. We further consider that both the prototype and the model are subjected to the gravitational fields $g$ and $\tilde{g}$, respectively; we quantity their ratio through the gravitational scale factor,
\begin{equation*}
\varsigma= \dfrac{\tilde{g}}{g}.
\end{equation*}
The density scale factor is identified by
\begin{equation*}
\gamma= \dfrac{\tilde{\rho}}{\rho}.
\end{equation*}
The friction coefficients are assumed to be equal in both systems$-$that is, $\lambda_{\mu}=\tilde{\mu}/\mu=1$. Notice that this assumption is made to simplify the following developments, but it is not stricly necessary.

Following the similarity statements determined by the $\pi$ terms, we can determine how, derived quantities are scaled in the model, upon the definition of the scaling factors $\lambda$, $\varsigma$, and $\gamma$. The scaling factor of a quantity `$f$' is identified as $\lambda_f=\tilde{f}/{f}$.
The scaling laws for the rigid-body motion require that
\begin{subequations}
	\begin{align}
	\dfrac{\tilde{J}_i}{J_i}=1, \quad \lambda_x = \lambda, \quad	\lambda_{\theta}= 1,\\
	\lambda_{\ddot{\theta}} = \dfrac{\varsigma}{\lambda},\quad
	\lambda_{\dot{\theta}} = \sqrt{\dfrac{\varsigma}{\lambda}}, \quad 
	\lambda_{\ddot{x}} = \varsigma,\quad
	\lambda_{\dot{x}} = \sqrt{\varsigma\lambda},\quad
	\lambda_t = \sqrt{\varsigma\lambda},\\
	\lambda_{\mu}=1,\quad
	\lambda_{J} = \gamma \lambda^5, \quad
	\lambda_{\mathcal{P}}= \gamma \varsigma \lambda,\quad
	\lambda_{\mathcal{I}}= \gamma \sqrt{\varsigma \lambda^3}.
\end{align}
\end{subequations}
We emphasize that, if a geometric scaling is imposed, only two parameters, $\lambda_{\mathcal{P}}$ and $\lambda_{\mathcal{I}}$, need to be specified \citep[see also ][]{199173}. The independent scaling factors are three, i.e., $\lambda$, $\varsigma$, and $\gamma$, but only two equations ($\lambda_{\mathcal{P}}= \gamma \varsigma \lambda$ and $\lambda_{\mathcal{I}}= \gamma \sqrt{\varsigma \lambda^3}$) need to be satisfied in defining the model. The system is thus over-determined, as the three scaling factors cannot independently satisfy the expressions of $\lambda_{\mathcal{P}}$ and $\lambda_{\mathcal{I}}$. This renders the derivation of scaling laws for blast actions non-trivial and challenging.
\section{Scaling laws}
\label{sec:scaling}
\noindent We assume that both the model and the prototype share the same gravitational field, so that $\varsigma=1$. Two parameters have to be selected: the geometric scaling, $\lambda$, and the density scaling factor, $\gamma$.\\
As previously discussed, both the duration time, the pressure peaks, and the impulse are functions of the stand-off distance, $R$, and the explosive quantity, $W$. In particular, $P_{ro}= \hat{P}_{ro}(Z)$ is function of the scaled distance, $Z$, while $t_{o}$ and $i_{ro}$ are functions of both $Z$ and $W$, i.e., $t_{o}= W^{1/3} \hat{t}_{ow}(Z)$ and $i_{ro}= W^{1/3} \hat{i}_{row}(Z)$. The above functions are presented in Appendix I.\\
By definition of the scaled distance $Z=R/W^{1/3}$, the explosive quantity scaling factor reads $\lambda_W = (\lambda/\lambda_Z)^3$. The scaling factors for $i_{ro}$, $t_o$, and $P_{ro}$ are hence:
\begin{equation}
	\lambda_{i_{ro}} = \dfrac{\hat{i}_{row}(\tilde{Z})}{\hat{i}_{row}\left({Z}\right)} \lambda_w^{\frac{1}{3}}, \quad 	
	\lambda_{t_{o}} = \dfrac{\hat{t}_{ow}(\tilde{Z})}{\hat{t}_{ow}\left({Z}\right)} \lambda_w^{\frac{1}{3}}, \quad
	\lambda_{P_{ro}} = \dfrac{\hat{P}_{ro}(\tilde{Z})}{\hat{P}_{ro}\left({Z}\right)}.
	\label{eq:similW_iro}
\end{equation}
The selection of $\lambda_Z-$that is, $\lambda_W$, as the geometric scaling for the stand-off distance is here imposed to be equal to $\lambda$. Therefore, the following identities must hold
\begin{subequations}
	\begin{align}
	\lambda_{i_{ro}} = \lambda_{\mathcal{I}}= \gamma \sqrt{\lambda^3},\quad
	\lambda_{t_{o}} = \lambda_{t}= \sqrt{\lambda},\quad
	\lambda_{P_{ro}} = \lambda_{\mathcal{P}}= \gamma \lambda.
	\label{eq:simil_iro}
\end{align}
\end{subequations}
It can be proven that no possible solution exists for $\lambda_Z$ such that all the similarity statements are verified simultaneously \cite[see ][]{199173}. Nevertheless, the system can be relaxed if one considers that the blast load is fast enough, compared to the characteristic time of the structure. In this case the blast load is considered as an impulsive load. This approximation is usually true in a wide range of applications, nevertheless it requires to be a posteriori verified.

For impulsive loads, $\pi_{23}$ vanishes (i.e., $\pi_{23}<<1$) and any difference between $\lambda_{t_o}$ and $\lambda_t$ is negligible as the time-history of the load is no more a main parameter (impulsive loading hypothesis). Therefore, from the three initial equations, only one needs to be verified$-$that is, $\lambda_{i_{ro}} = \lambda_{\mathcal{I}}$. Due to the high non-linearity of the function $\hat{i}_{row}(Z)$ (see Appendix I), an analytical solution is not obvious. The scaling factor for the scaled distance can be found by solving the following non-linear equation:
\begin{equation}
\text{find } \lambda_Z \text{ such that } \qquad \dfrac{1}{\lambda_Z}\dfrac{\hat{i}_{row}(Z\lambda_Z)}{\hat{i}_{row}(Z)}= \gamma \sqrt{\lambda}.
\label{eq:simil_z}
\end{equation}
It is worth mentioning that there exists a particular case of scaling laws for which $\lambda_Z= \hat{Z}/Z=1$. This is the Hopkinson-Cranz similarity law \cite{hopki,cranz}. In this case, the scaling factor for the explosive quantity is $\lambda_W=\lambda^3$. Nevertheless, differently from the scaling proposed herein, the Hopkinson-Cranz similarity law prescribes the density scaling factor: $\gamma = 1/\sqrt{\lambda}$, which is an important and impractical restriction.

The general scaling law, Eq. (\ref{eq:simil_z}), does not have restrictions neither on the geometric scaling, $\lambda$, nor on the density (or mass) scaling, $\gamma$. For the particular (Hopkinson-Cranz) case, the mass of the model, instead, is directly identified as $\tilde{m}=m \lambda^{5/2}$. In order to obtain such values, either the model material should have much higher density than that of the prototype material or masses should be added for assuring the proper equivalent density, by respecting the mass moment inertia similarity. This is hard to achieve in practice. Moreover, in this second scenario, the scaled distance in the model would equal the scaled distance in the prototype, which means that in both systems the pressure $P_{ro}$ would have to be the same. This usually represents a disadvantageous scaling, as one of the main objectives of conducting in-scale experimental tests is the reduction of the intensity of the explosive load.\\
The proposed scaling laws and the case of Hopkinson-Cranz similarity are summarized in Table \ref{tab:scaling_general} and \ref{tab:scaling_particular}, respectively.\\

It is worth noticing that, relying on similar developments, scaling laws accounting also for the material response (elastic and inelastic) can be obtained, as detailed in \citet{masithesis}, among others.

\begin{table}
\caption{Relations for model and prototype variables, for general case.}
\label{tab:scaling_general}
\centering
\small
\renewcommand{\arraystretch}{1.25}
\begin{adjustbox}{max width=1\textwidth}
\begin{tabular}{l c l c}
\hline\hline
			Variable$\quad$ & Scaling factor& Variable$\quad$ & Scaling factor\\
			\hline
			Length, $l$ & $\lambda$ & Angle, $\theta$ & 1 \\
			Material density, $\rho$ & $\gamma$ & Angular velocity, $\dot{\theta}$ & $\lambda^{-1/2}$\\
			Linear displacement, $x$ & $\lambda$ & Angular acceleration, $\ddot{\theta}\quad $ & 1\\
			Linear velocity, $\dot{x}$& $\lambda^{1/2}$ & Time, $t$ & $\lambda^{1/2}$ \\
			Linear acceleration, $\ddot{x}\qquad $  & $\lambda^{-1}$& Mass, $m$ & $\gamma \lambda^3$\\
			Blast impulse & $\gamma \lambda^{3/2}$  & Mass moment of inertia, $J$ & $\gamma \lambda^5$ \\
			TNT equivalent, $W$ & Eq. (\ref{eq:simil_z})\\
\hline\hline
\end{tabular}
\end{adjustbox}
\normalsize
\end{table}

\begin{table}
\caption{Relations for model and prototype variables, respecting Hopkinson-Cranz similarity.}
\label{tab:scaling_particular}
\centering
\small
\renewcommand{\arraystretch}{1.25}
\begin{adjustbox}{max width=1\textwidth}
\begin{tabular}{l c l c}
\hline\hline
			Variable$\quad$ & Scaling factor& Variable$\quad$ & Scaling factor\\
			\hline
			Length, $l$ & $\lambda$ & Angle, $\theta$ & 1 \\
			Material density, $\rho$ & $\lambda^{-1/2}$ & Angular velocity, $\dot{\theta}$ & $\lambda^{-1/2}$\\
			Linear displacement, $x$ & $\lambda$ & Angular acceleration, $\ddot{\theta}\quad $ & 1\\
			Linear velocity, $\dot{x}$& $\lambda^{1/2}$ & Time, $t$ & $\lambda^{1/2}$ \\
			Linear acceleration, $\ddot{x}\qquad $  & $\lambda^{-1}$& Mass, $m$ & $\lambda^{5/2}$\\
			Blast impulse & $ \lambda$ & Mass moment of inertia, $J$ & $\lambda^{9/2}$ \\
			TNT equivalent, $W$ & $\lambda^3$\\
\hline\hline
\end{tabular}
\end{adjustbox}
\normalsize
\end{table}

\section{Validation of the scaling laws}
\label{sec:sim_validation}
\noindent Once the scaling laws for the structural dynamic response and blast loading have been derived, we investigate their validity for masonry structures displaying monolithic behavior.\\
First, we compare the prototype and model responses obtained by numerical integration of the non-linear equations of rocking motion \citep[see ][]{rockingmasi}. Then, the validation is performed relying on three-dimensional Finite Element (FE) simulations. 

The target consists of a rectangular (rigid) block with uniformly distributed mass $m$. The dimensions of the block are $2 b \times 2 h \times 2 w$ and the radial distance from the rocking pivot point $O$ to the center of gravity is $r=h/\cos \alpha$, with $\alpha$ being the slenderness angle. The equations of motion for a rocking response mechanism are
\begin{equation}
	\mathcal{J}_o \ddot{\theta}+m g r \sin \left[\alpha\;\text{sgn}\left(\theta\right)-\theta\right]= S r P_{r} (t)  \cos \left[\alpha\;\text{sgn}\left(\theta\right)-\theta\right],
\label{eq:eom} 
\end{equation}
where $\mathcal{J}_o=\frac{4}{3} m r^2$ is the moment of inertia with respect to the pivot point and $\theta=\theta(t)$ is the inclination angle \citep[see ][]{rockingmasi,lourenco1,makris,dejong,elia,masi2020resistance}.\\
In order to validate the general scaling law, summarized in Table \ref{tab:scaling_general}, which is different from the Hopkinson-Cranz one (Table \ref{tab:scaling_particular}), we consider several case studies. First, the similarity between prototype and models is tested considering similar materials, i.e., $\gamma=1$. Then, the case of different materials, with different density (i.e., $\gamma \neq 1$), is investigated.\\

Finally, the validation is performed through three-dimensional numerical Finite Element (FE) simulations. The applicability of the proposed scaling law is thus investigated under the combined effects of sliding, rocking, and uplift (flight mode). We consider Coulomb friction at the interface of the block with the rigid base, with an angle of friction equal to $35^{\circ}$, which is common for many geomaterials (concrete, marble, stone etc.). Blast loads are applied relying on the best-fit interpolations in Appendix I. ABAQUS commercial software is used for the computations \cite{abaqus}. A hard contact formulation is used, i.e., no penetration is allowed at the contact of the rocking block with the base. The rigid base is fixed and the rigid block is free to translate and rotate along all directions (see Figure \ref{f:rectangular_block_model}).

\subsection{Equations of motions}
\label{subsection:eom}
\subsubsection{Same material modeling, $\gamma=1$}
\label{par:similar}
\begin{figure*}[h]
	\begin{center}
		\includegraphics[width=0.65\textwidth]{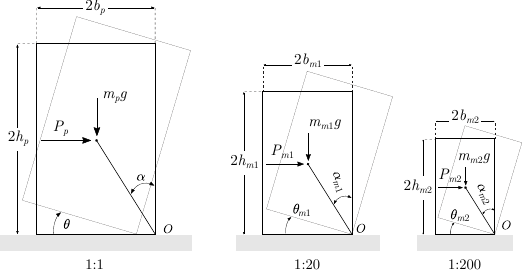}
		\caption{Prototype system (left) and models with geometric scaling $\lambda=1/20$ (center) and $\lambda=1/200$ (right).}
		\label{fig:scaling}
	\end{center}
\end{figure*}

\noindent We consider a prototype block, with density $\rho_p=2000$ kg m$^{-1/3}$, height $2h_p=10$ m, slenderness angle $\alpha_p=15^{\circ}$, and arbitrary depth $2w_p$, subjected to the loading of the detonation of a given explosive quantity, $W_p$, at a stand-off distance $R_p =2$ m. Two different geometric scales are considered; $\lambda=1/20$ and $\lambda=1/200$, see Figure \ref{fig:scaling}. Table \ref{tab:scaling_parameters} displays the geometry parameters for the prototype and the two models. For both models, we assume a unit density scaling factor, i.e., $\gamma=1$. Three different quantities of TNT equivalent are considered: (a) $W_{pa} = 50$ kg, (b) $W_{pb} = 100$ kg, and (c) $W_{pc} = 79.8$ kg. For case (a), according to the analytical developments in \citet{rockingmasi}, and by numerical integration of Eq. (\ref{eq:eom}), the prototype block rocks, without overturning. Overturning is instead expected for case (b). Case (c), $W_{pc} = 79.8$ kg, represents the maximum critical quantity of explosive at the stand-off distance $R_p=2$ m for which toppling does not happen.

By numerically solving Eq. (\ref{eq:simil_iro}) for $\lambda_Z$ and for each loading case scenario (a-c), we obtain the quantities of explosive (and scaled distances) for the models, as shown in Table \ref{tab:scaling_TNT}.

\begin{table}
\caption{Geometry parameters for the prototype and the two models, shown in Figure \ref{fig:scaling}.}
\label{tab:scaling_parameters}
\centering
\small
\renewcommand{\arraystretch}{1.25}
\begin{adjustbox}{max width=1\textwidth}
\begin{tabular}{l l c c c l c c c l l c c}
\hline\hline
			Prototype&& &&$\quad$&Model 1&& $\quad$& &Model 2 &&\\
			\cline{1-4} \cline{6-8} \cline{10-12}
			Height& $2h_p$ &(m) & 10 & &$h_{m1}$ &(cm)& 50 & &$h_{m2}$&(cm) & 5\\
			Width& $2b_p$ &(m) & 2.68&&$2b_{m1}$ &(cm)&13.40&&$2b_{m2}$&(cm)&1.34\\
			Slenderness angle& $\alpha_p$ &($^{\circ}$) & 15 & &$\alpha_{m1}$ &($^{\circ}$)& 15 & &$\alpha_{m2}$&($^{\circ}$) & 15\\
			Stand-off distance& $R_p$ &(m) & 2& & $R_{m1}$ &(cm)& 10 &  & $R_{m2}$&(cm)& 1\\
\hline\hline
\end{tabular}
\end{adjustbox}
\normalsize
\end{table}

\begin{table}[hbt]
\caption{Blast loading characteristics for the prototype and the two models, shown in Figure \ref{fig:scaling}.}
\label{tab:scaling_TNT}
\centering
\small
\renewcommand{\arraystretch}{1.25}
\begin{adjustbox}{max width=1\textwidth}
\begin{tabular}{l c l c l c l c l c l c c}
\hline\hline
			Prototype& &&$\quad$&Model 1&&& $\quad$&Model 2 &&\\
			\cline{1-3} \cline{5-7} \cline{9-11}
			TNT, $W_{pa}$ & 50 &kg & &$W_{m1}$ & 1.0 &mg & &$W_{m2}$ & 0.233& $\mu$g\\
			Scaled distance, $Z_{pa}$ & 0.54 &m$\sqrt[3]{\text{kg}}$ & &$Z_{m1}$ & 0.99 &m$\sqrt[3]{\text{kg}}$ & &$Z_{m2}$ & 1.62& m$\sqrt[3]{\text{kg}}$\\
			\hline
			TNT, $W_{pb}$ & 100 &kg & &$W_{m1}$ & 2.06 &mg & &$W_{m2}$ & 0.492& $\mu$g\\
			Scaled distance, $Z_{pb}$ & 0.43 &m$\sqrt[3]{\text{kg}}$ & &$Z_{m1}$ & 0.78 &m$\sqrt[3]{\text{kg}}$ & &$Z_{m2}$ & 1.26& m$\sqrt[3]{\text{kg}}$\\
			\hline
			TNT, $W_{pc}$ & 79.8 &kg & &$W_{m1}$ & 1.635 &mg & &$W_{m2}$ & 0.386& $\mu$g\\
			Scaled distance, $Z_{pc}$ & 0.46 &m$\sqrt[3]{\text{kg}}$ & &$Z_{m1}$ & 0.85 &m$\sqrt[3]{\text{kg}}$ & &$Z_{m2}$ & 1.37& m$\sqrt[3]{\text{kg}}$\\
\hline\hline
\end{tabular}
\end{adjustbox}
\normalsize
\end{table}

\begin{figure*}[ht]
\centering
\begin{subfigure}[b]{0.3\textwidth}
  \centering
  \includegraphics[width=\linewidth]{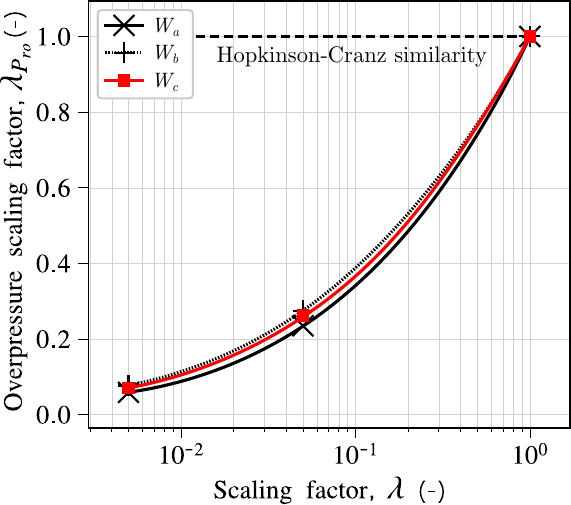}
	\caption{\footnotesize overpressure}
	\label{f:scaling_adimensional_a}
\end{subfigure}\hspace{0.2cm}
\begin{subfigure}[b]{0.3\textwidth}
  \centering
  \includegraphics[width=\linewidth]{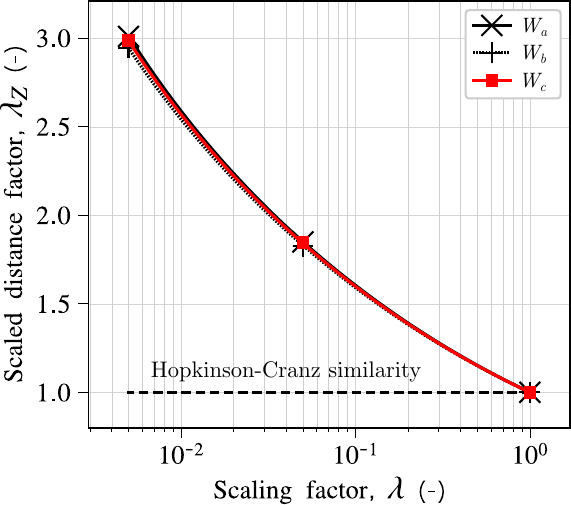}
	\caption{\footnotesize scaled distance}
	\label{f:scaling_adimensional_b}
\end{subfigure} \hspace{0.2cm}
\begin{subfigure}[b]{0.3\textwidth}
  \centering
  \includegraphics[width=\linewidth]{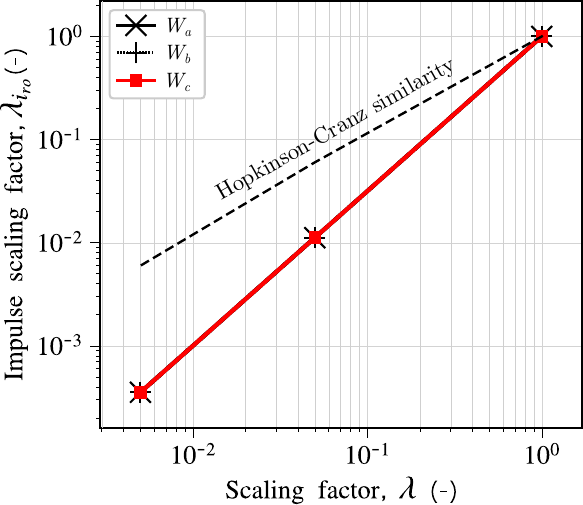}
	\caption{\footnotesize impulse}
	\label{f:scaling_adimensional_c}
\end{subfigure}
\caption{Scaling factors of overpressure (a), scaled distance (b), and impulse (c), in function of the geometric scaling factor, $\lambda$ (see Table \ref{tab:scaling_parameters}). The scaling law allows higher overpressure and impulse reduction than the Hopkinson-Cranz similarity law (see Table \ref{tab:scaling_particular}).}
	\label{f:scaling_adimensional}
\end{figure*}

We can notice that for scaling factors $\lambda<1$, the calculated scaled distance of the model is higher than that one of the prototype. This is a favorable feature of the proposed scaling, since the intensity of the blast in the model is smaller than that in the prototype. The overpressure peak and impulse have thus smaller intensities in the model. This is not the case for Hopkinson-Cranz scaling law. Figure \ref{f:scaling_adimensional} displays the dependency of the scaling factors on overpressure (\ref{f:scaling_adimensional_a}), scaled distance (\ref{f:scaling_adimensional_b}), and impulse (\ref{f:scaling_adimensional_c}) with respect to the geometric scaling. The Hopkinson-Cranz scaling is shown by dashed lines.\\
For scaling factor $\lambda=1/200$, the overpressure peak is only $5\div 8\%$ the value for the prototype. The impulse in the model is found to be only $0.3\%$ the impulse in the prototype. Instead, if the Hopkison-Cranz law is used, the model overpressure peak equals the prototype value and the impulse is as large as $5\%$ the value in the prototype (one order of magnitude higher).\\

\begin{figure*}[t]
\centering
\begin{subfigure}[b]{0.3\textwidth}
  \centering
  \includegraphics[width=\linewidth]{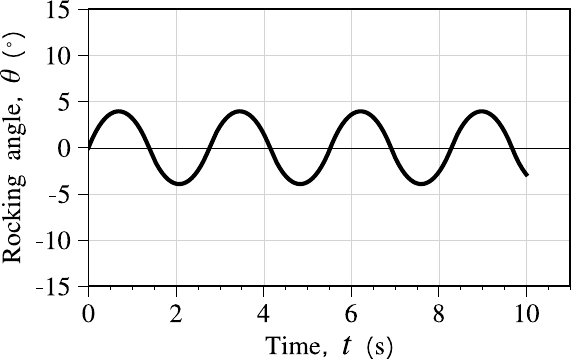}
	\caption{\footnotesize $W_{pa}$, rocking angle}
\end{subfigure}\hspace{0.1cm}
\begin{subfigure}[b]{0.3\textwidth}
  \centering
  \includegraphics[width=\linewidth]{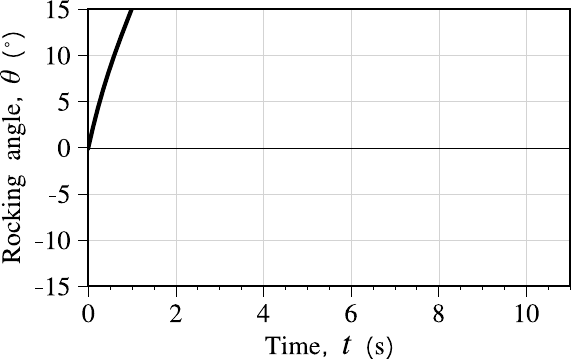}
	\caption{\footnotesize $W_{pb}$, rocking angle}
\end{subfigure} \hspace{0.1cm}
\begin{subfigure}[b]{0.3\textwidth}
  \centering
  \includegraphics[width=\linewidth]{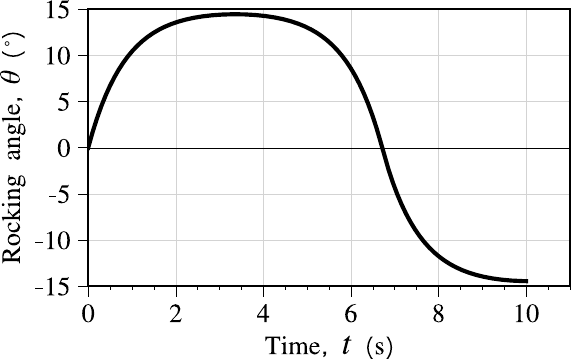}
	\caption{\footnotesize $W_{pc}$, rocking angle}
\end{subfigure}\\ \vspace{0.2cm}
\begin{subfigure}[b]{0.3\textwidth}
  \centering
  \includegraphics[width=\linewidth]{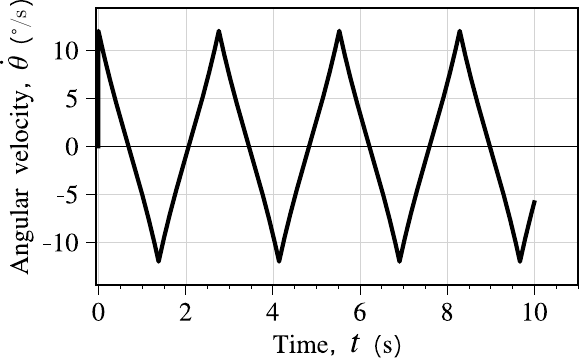}
	\caption{\footnotesize $W_{pa}$, angular velocity}
\end{subfigure}\hspace{0.1cm}
\begin{subfigure}[b]{0.3\textwidth}
  \centering
  \includegraphics[width=\linewidth]{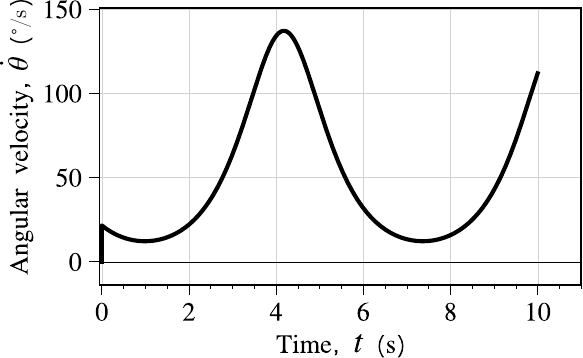}
	\caption{\footnotesize $W_{pb}$, angular velocity}
\end{subfigure} \hspace{0.1cm}
\begin{subfigure}[b]{0.3\textwidth}
  \centering
  \includegraphics[width=\linewidth]{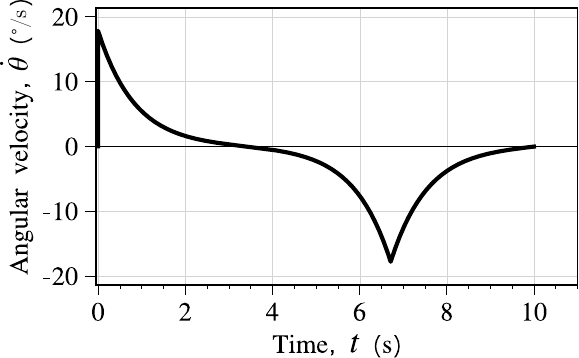}
	\caption{\footnotesize $W_{pc}$, angular velocity}
\end{subfigure}
\caption{Response of the prototype for $W_{pa}=50$ kg, $W_{pb}=100$ kg, and $W_{pc}=79.8$ kg, in terms of (a-c) the rocking angle $\theta$  and (d-f) the angular velocity $\dot{\theta}$.}
	\label{f:prototype_results}
\end{figure*}

\begin{figure*}[t]
\centering
\begin{subfigure}[b]{0.3\textwidth}
  \centering
  \includegraphics[width=\linewidth]{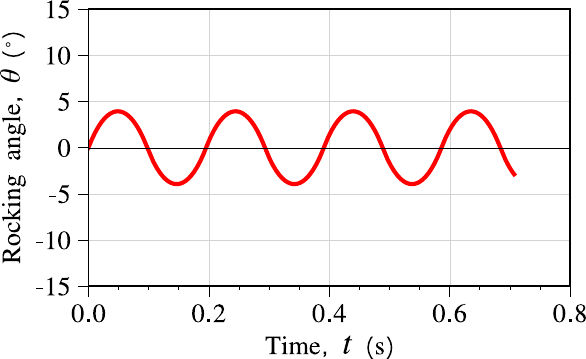}
	\caption{\footnotesize $W_{pa}$, rocking angle}
\end{subfigure}\hspace{0.1cm}
\begin{subfigure}[b]{0.3\textwidth}
  \centering
  \includegraphics[width=\linewidth]{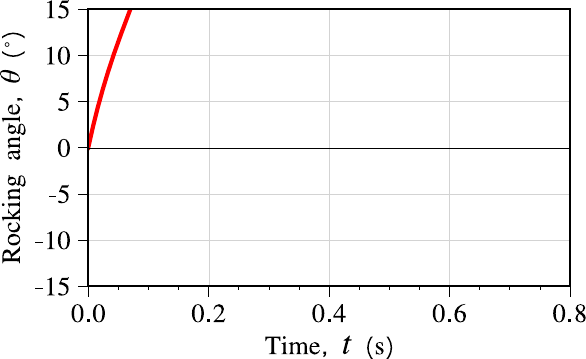}
	\caption{\footnotesize $W_{pb}$, rocking angle}
\end{subfigure} \hspace{0.1cm}
\begin{subfigure}[b]{0.3\textwidth}
  \centering
  \includegraphics[width=\linewidth]{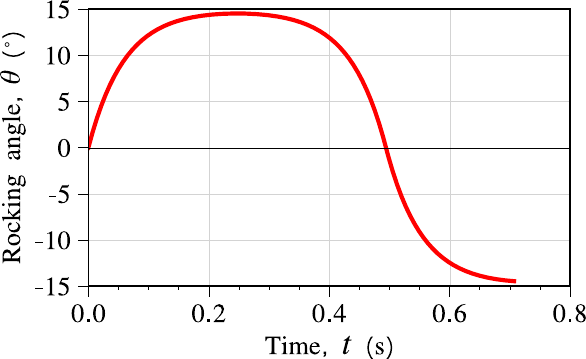}
	\caption{\footnotesize $W_{pc}$, rocking angle}
\end{subfigure}\\ \vspace{0.2cm}
\begin{subfigure}[b]{0.3\textwidth}
  \centering
  \includegraphics[width=\linewidth]{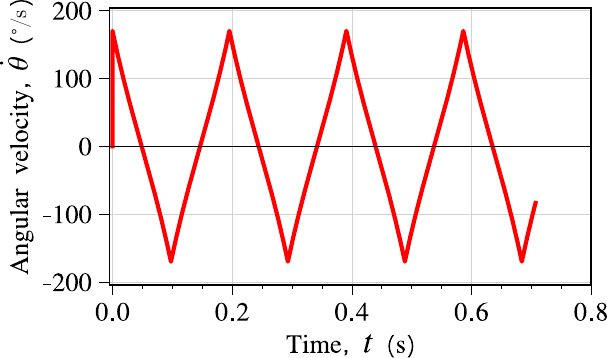}
	\caption{\footnotesize $W_{pa}$, angular velocity}
\end{subfigure}\hspace{0.1cm}
\begin{subfigure}[b]{0.3\textwidth}
  \centering
  \includegraphics[width=\linewidth]{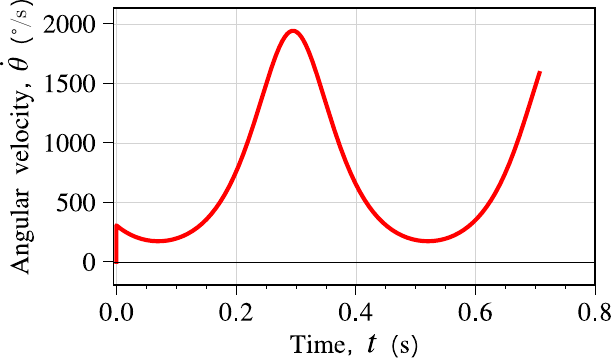}
	\caption{\footnotesize $W_{pb}$, angular velocity}
\end{subfigure} \hspace{0.1cm}
\begin{subfigure}[b]{0.3\textwidth}
  \centering
  \includegraphics[width=\linewidth]{fig_6d.pdf}
	\caption{\footnotesize $W_{pc}$, angular velocity}
\end{subfigure}
\caption{Response of the model for prototype explosive charges $W_{pa}=50$ kg, $W_{pb}=100$ kg, and $W_{pc}=79.8$ kg, in terms of (a-c) the rocking angle $\theta$  and (d-f) the angular velocity $\dot{\theta}$.}
	\label{f:model_results}
\end{figure*}

\begin{figure*}[ht]
\centering
\begin{subfigure}[h]{0.3\textwidth}
  \centering
  \includegraphics[width=\linewidth]{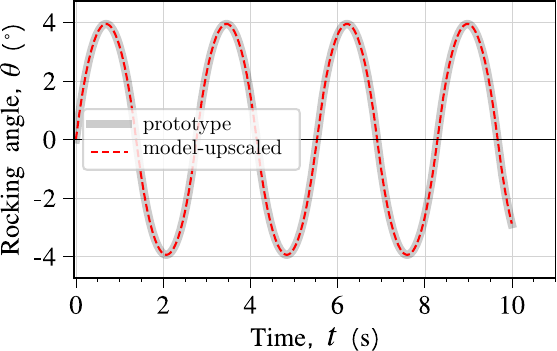}
	\caption{\footnotesize rocking angle, $W_{pa}=50$ kg}
\end{subfigure}\hspace{0.2cm}
\begin{subfigure}[h]{0.3\textwidth}
  \centering
  \includegraphics[width=\linewidth]{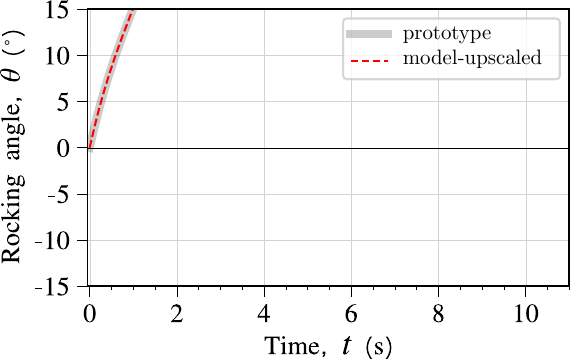}
	\caption{\footnotesize rocking angle, $W_{pb}=100$ kg}
\end{subfigure}\hspace{0.2cm}
\begin{subfigure}[h]{0.3\textwidth}
  \centering
  \includegraphics[width=\linewidth]{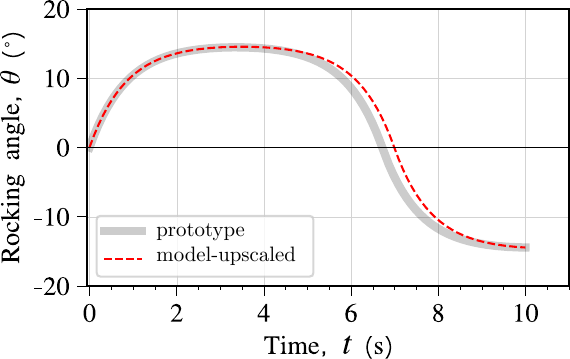}
	\caption{\footnotesize rocking angle, $W_{pc}=79.8$ kg}
\end{subfigure} \\
\vspace{0.2cm}
\begin{subfigure}[h]{0.3\textwidth}
  \centering
  \includegraphics[width=\linewidth]{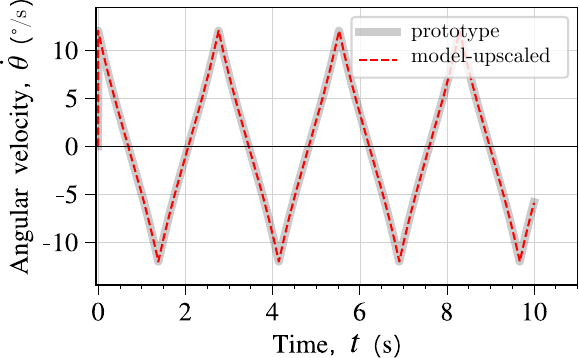}
	\caption{\footnotesize angular velocity, $W_{pa}=50$ kg}
\end{subfigure}\hspace{0.2cm}
\begin{subfigure}[h]{0.3\textwidth}
  \centering
  \includegraphics[width=\linewidth]{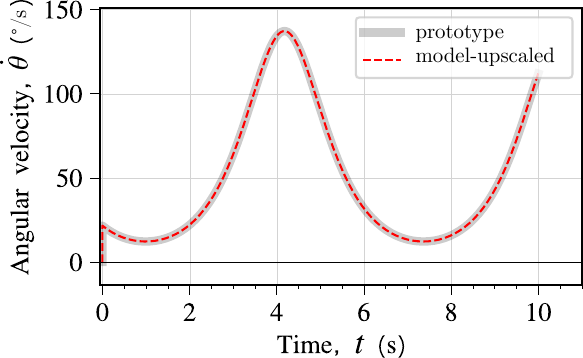}
	\caption{\footnotesize angular velocity, $W_{pb}=100$ kg}
\end{subfigure}\hspace{0.2cm}
\begin{subfigure}[h]{0.3\textwidth}
  \centering
  \includegraphics[width=\linewidth]{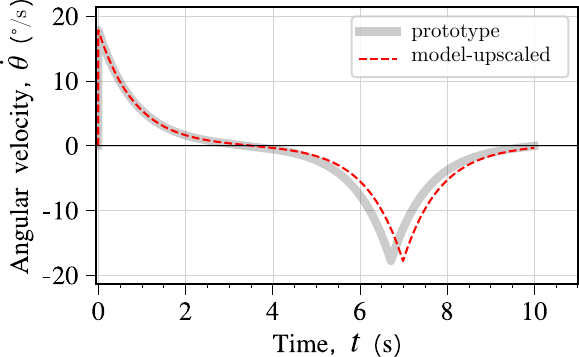}
	\caption{\footnotesize angular velocity, $W_{pc}=79.8$ kg}
\end{subfigure} \\
\vspace{0.2cm}
\begin{subfigure}[h]{0.3\textwidth}
  \centering
  \includegraphics[width=\linewidth]{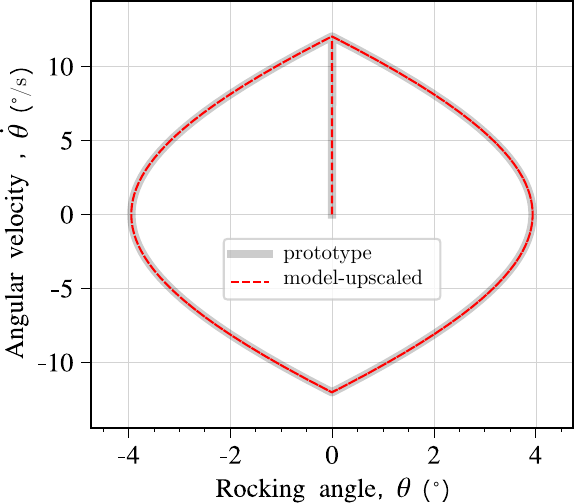}
	\caption{\footnotesize phase portrait, $W_{pa}=50$ kg}
\end{subfigure}\hspace{0.2cm}
\begin{subfigure}[h]{0.3\textwidth}
  \centering
  \includegraphics[width=\linewidth]{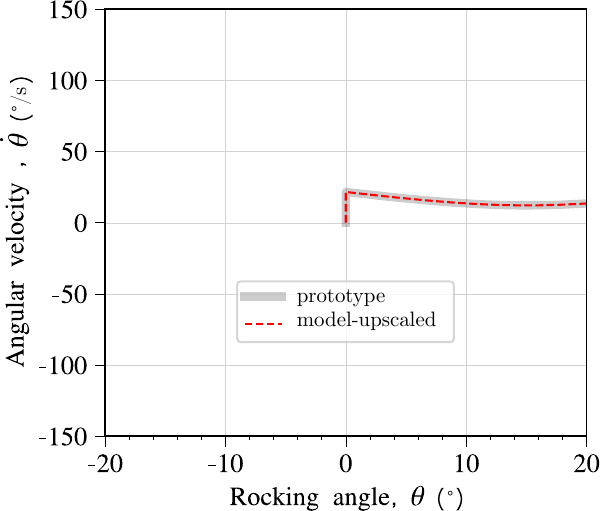}
	\caption{\footnotesize phase portrait, $W_{pb}=100$ kg}
\end{subfigure}\hspace{0.2cm}
\begin{subfigure}[h]{0.3\textwidth}
  \centering
  \includegraphics[width=\linewidth]{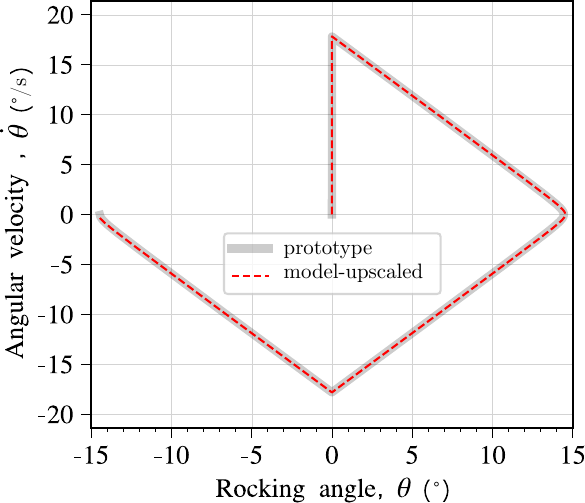}
	\caption{\footnotesize phase portrait, $W_{pc}=79.8$ kg}
\end{subfigure}
\caption{Comparison between the prototype response and the model response ($\lambda=1/200$). The model is upscaled, i.e., all quantities are multiplied by the inverse of the scaling factor (cf. Table \ref{tab:scaling_parameters}). Columns represent the three different quantities of explosive: $W_{pa}$ (left), $W_{pb}$ (center), and $W_{pc}$ (right).}
	\label{f:scaling_upscaled}
\end{figure*}

From numerical integration of (non-linear) Eq. (\ref{eq:eom}), we compare the response of the prototype with the response of model 2, which has $\lambda=1/200$. Figure \ref{f:prototype_results} displays the evolution of the rocking angle and the angular velocity for the prototype system and the three quantities of explosive. Figure \ref{f:model_results} shows the response of the scaled model. The scaled model agrees with the prototype in terms of the final state of the block. For $W_{pa}$, both systems rock (without toppling); for $W_{pb}$, both systems undergo overturning; while $W_{pc}$ and its scaled counterpart represent the critical explosive quantity of both systems.\\
Figure \ref{f:scaling_upscaled} compares the prototype response with that of the model, upscaled, i.e., all quantities are multiplied by the inverse of the scaling factor (cf. Table \ref{tab:scaling_parameters}). The curves of the systems coincide which confirms the derived scaling laws. For the critical explosive quantity, $W_{pc}$, a negligible offset between the model and the prototype exists. This is a special, critical case, as it refers to the critical explosive quantity to avoid overturning. In this case, numerical errors may become important and can determine overturning or not. Nevertheless, the scaling law is found to correctly capture the dynamics of the prototype when the phase space is examined.

\subsubsection{Different material modeling, $gamma \neq 1$}
\label{par:dissimilar}
\begin{figure*}[ht]
\centering
\begin{subfigure}[b]{0.3\textwidth}
  \centering
  \includegraphics[height=0.85\linewidth]{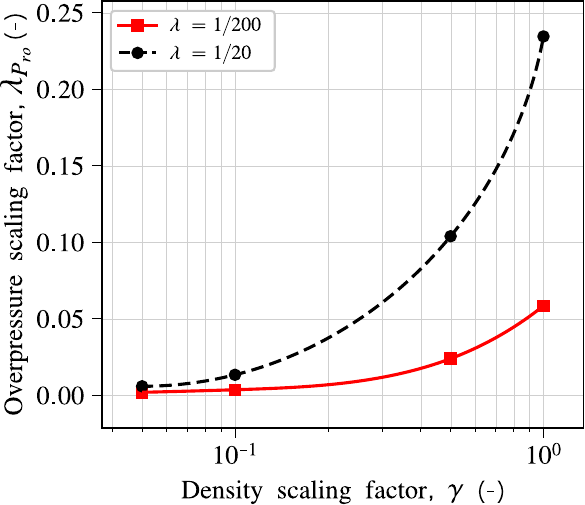}
	\caption{\footnotesize scaling factor}
\end{subfigure}\hspace{0.2cm}
\begin{subfigure}[b]{0.3\textwidth}
  \centering
  \includegraphics[height=0.85\linewidth]{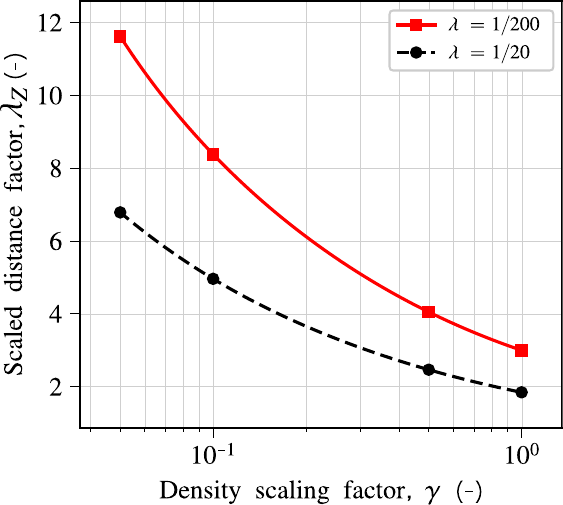}
	\caption{\footnotesize scaled distance}
\end{subfigure}
\begin{subfigure}[b]{0.3\textwidth}
  \centering
  \includegraphics[height=0.85\linewidth]{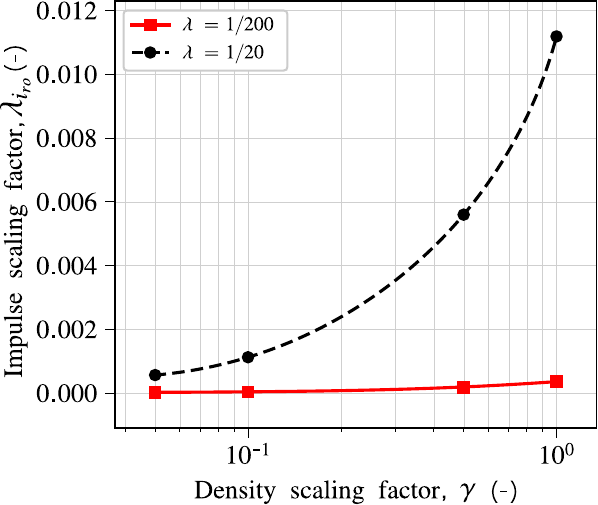}
	\caption{\footnotesize impulse}
\end{subfigure}
\caption{Scaling factors of overpressure (a), scaled distance (b), and impulse (c), in function of the density scaling factor, $\gamma$ (see Table \ref{tab:scaling_parameters}), and for two geometric scaling factors and prototype explosive quantity 50 kg.}
	\label{f:scaling_dissimilar_tend}
\end{figure*}

\noindent We consider here the scaling laws in the frame of different material modeling. This may be advantageous in experimental tests, as it allows to reduce the intensity of the blast load (cf. Table \ref{tab:scaling_general}). More specifically, a material model with density lower than that of the prototype material allows reducing both the model pressure peak and impulse. We emphasize that this holds true only for the proposed scaling laws (cf. Hopkinson-Cranz, Table \ref{tab:scaling_particular}). In Figure \ref{f:scaling_dissimilar_tend} we show how the overpressure, scaled distance, and impulse scaling factors vary for a given geometric scaling. For instance, assuming $\lambda=1/200$, the model overpressure peak equals $5.8$\%  the prototype peak, with a density scaling factor $\gamma=1$. For the same scaling, but with a dissimilar material such that $\gamma=0.05$, the model overpressure peak reduces to $0.18$\%.

\begin{figure*}[b]
\centering
\begin{subfigure}[h]{0.3\textwidth}
  \centering
  \includegraphics[width=\linewidth]{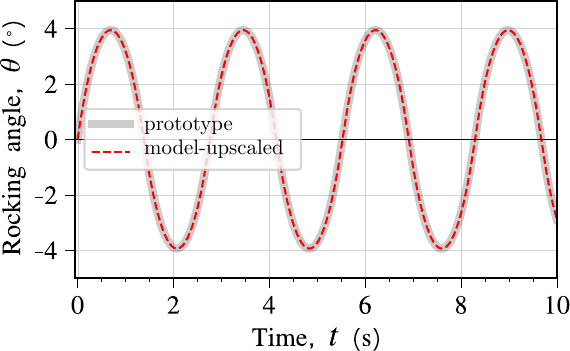}
	\caption{\footnotesize rocking angle, $W_{pa}=50$ kg}
\end{subfigure}\hspace{0.2cm}
\begin{subfigure}[h]{0.3\textwidth}
  \centering
  \includegraphics[width=\linewidth]{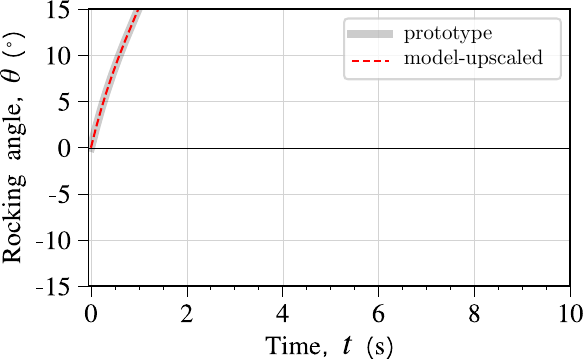}
	\caption{\footnotesize rocking angle, $W_{pb}=100$ kg}
\end{subfigure}\hspace{0.2cm}
\begin{subfigure}[h]{0.3\textwidth}
  \centering
  \includegraphics[width=\linewidth]{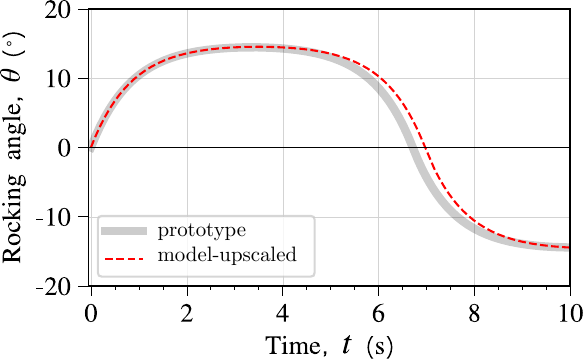}
	\caption{\footnotesize rocking angle, $W_{pc}=79.8$ kg}
\end{subfigure} \\
\vspace{0.2cm}
\begin{subfigure}[h]{0.3\textwidth}
  \centering
  \includegraphics[width=\linewidth]{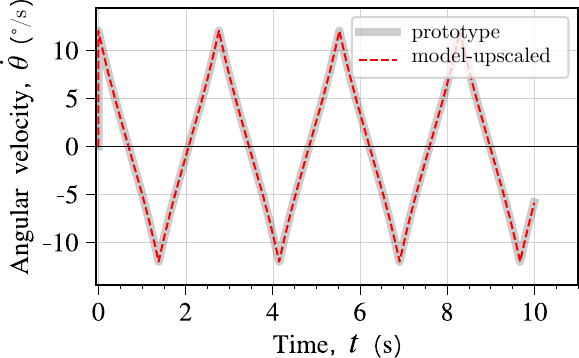}
	\caption{\footnotesize angular velocity, $W_{pa}=50$ kg}
\end{subfigure}\hspace{0.2cm}
\begin{subfigure}[h]{0.3\textwidth}
  \centering
  \includegraphics[width=\linewidth]{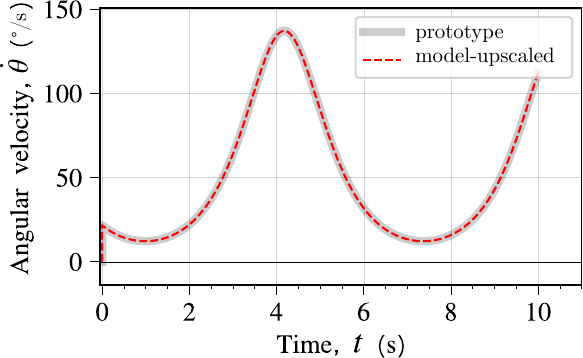}
	\caption{\footnotesize angular velocity, $W_{pb}=100$ kg}
\end{subfigure}\hspace{0.2cm}
\begin{subfigure}[h]{0.3\textwidth}
  \centering
  \includegraphics[width=\linewidth]{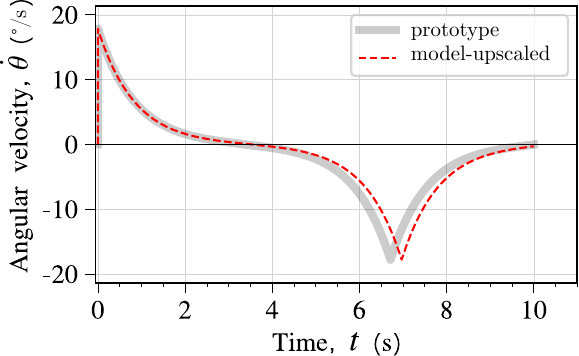}
	\caption{\footnotesize angular velocity, $W_{pc}=79.8$ kg}
\end{subfigure}
\caption{Comparison between the prototype response and the model response ($\lambda=1/200$) and dissimilar materials ($\gamma=0.07$). The model is upscaled, i.e., all quantities are multiplied by the inverse of the scaling factor (cf. Table \ref{tab:scaling_parameters}). Columns represent the three different quantities of explosive: $W_{pa}$ (left), $W_{pb}$ (center), and $W_{pc}$ (right).}
	\label{f:scaling_dissimilar}
\end{figure*}

As previously done for the similar material modeling, we validate, through the numerical integration of the equation of motion (\ref{eq:eom}), the scaling laws for materials with different densities. A light material as balsa (average density equal to $140$ kg m$^{-1/3}$) is considered for the model system. The density scaling factor is thus $\gamma=0.07$. In Figure \ref{f:scaling_dissimilar} the response of the prototype is compared with the up-scaled results of the model. The comparison gives the same results with the case of same materials between model and prototype.

\subsection{FEM simulations}
\label{subsection:fem}
\begin{figure}[ht]
\centering
  \includegraphics[width=0.2\linewidth]{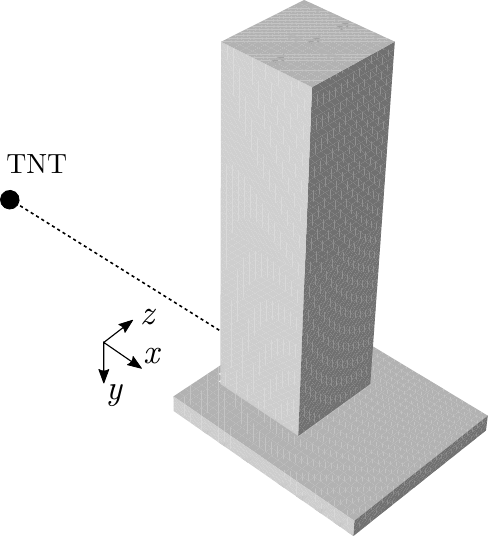}
\caption{Geometric and mesh discretization of the rigid block subjected to a TNT explosion.}
	\label{f:rectangular_block_model}
\end{figure}

\noindent After validation of the proposed scaling laws for rocking response mechanisms, we consider the more realistic scenario of a masonry structure undergoing rocking, sliding, and up-lifting. This is investigated through three-dimensional FE simulations of a rigid block standing on a rigid base, see Figure \ref{f:rectangular_block_model}. First, slender blocks, with geometric dimensions as in Table \ref{tab:scaling_parameters}, are considered. Then, non-slender blocks are also investigated.\\

\begin{figure*}[ht]
  \centering
  \begin{subfigure}[h]{0.5\textheight}
  \centering
  \includegraphics[width=\linewidth]{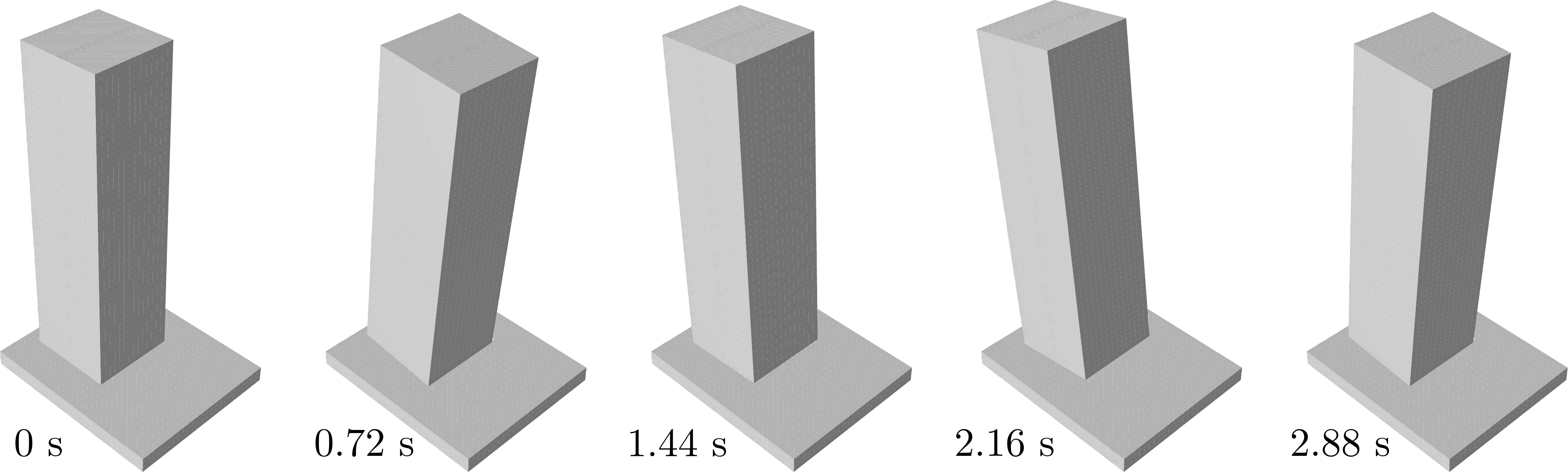}
	\caption{\footnotesize $W_{pa}=50$ kg}
\end{subfigure}\\
\vspace{0.2cm}
\begin{subfigure}[h]{0.5\textheight}
  \centering
  \includegraphics[width=0.8\linewidth]{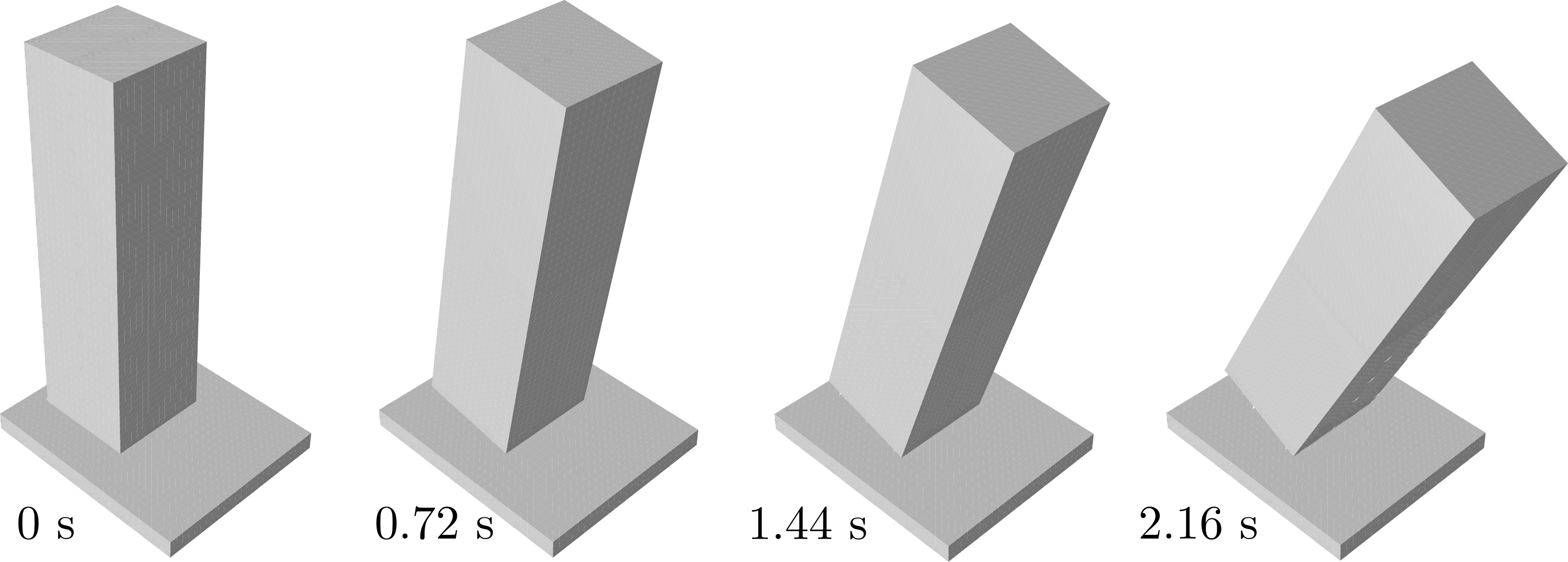}
	\caption{\footnotesize $W_{pb}=100$ kg}
\end{subfigure} 
\caption{Dynamic evolution of a rigid block subjected to (a) $W_{pa}=50$ kg and (b) $W_{pb}=100$ kg, according to the FE simulations.}
	\label{f:block_FE_sim}
\end{figure*}

In the case of slender blocks, two quantities of explosive are considered: $W_{pa}=50$ kg and $W_{pb}=100$ kg (prototype values). The evolution of the rigid-body motion is displayed in Figure \ref{f:block_FE_sim}. It is worth noticing that rocking prevails over the other mechanisms.\\
In the case of non-slender blocks, one quantity of explosive is only considered, i.e., 2500 kg. Even if such a large quantity of TNT may cause damage to the target (instead of a rigid-body response), this choice is only made for validating the scaling laws for mechanisms where sliding is dominant.

\subsubsection{Slender blocks}
\label{par:slenderblocks}
\noindent By applying the scaling laws for a geometric scaling factor $\lambda=1/200$, two different models are investigated. We consider a model made of the same material with the prototype. The block dimensions are those of Table \ref{tab:scaling_parameters} (Prototype and Model 2). \\

\begin{figure*}[h]
\centering
\begin{subfigure}[h]{0.3\textwidth}
  \centering
  \includegraphics[width=\linewidth]{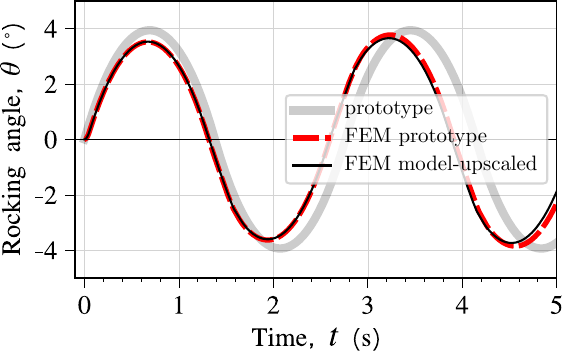}
	\caption{\footnotesize rocking angle, $W_{pa}=50$ kg}
\end{subfigure}\hspace{0.2cm}
\begin{subfigure}[h]{0.3\textwidth}
  \centering
  \includegraphics[width=\linewidth]{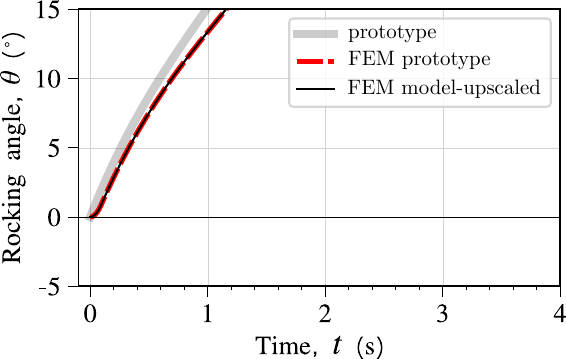}
	\caption{\footnotesize rocking angle, $W_{pb}=100$ kg}
\end{subfigure} \\
\vspace{0.2cm}
\begin{subfigure}[h]{0.3\textwidth}
  \centering
  \includegraphics[width=\linewidth]{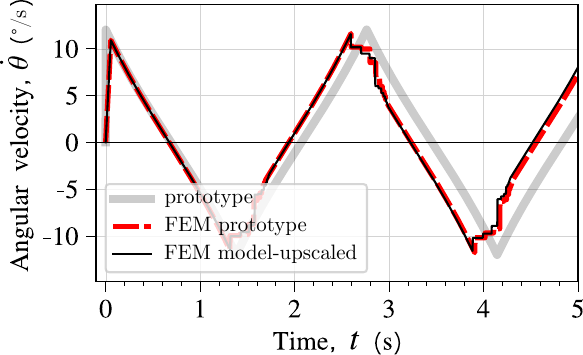}
	\caption{\footnotesize angular velocity, $W_{pa}=50$ kg}
\end{subfigure}\hspace{0.2cm}
\begin{subfigure}[h]{0.3\textwidth}
  \centering
  \includegraphics[width=\linewidth]{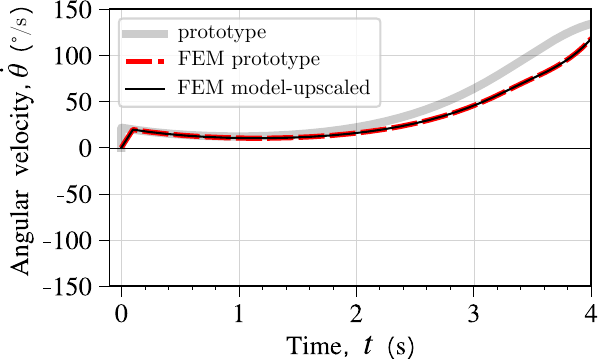}
	\caption{\footnotesize angular velocity, $W_{pb}=100$ kg}
\end{subfigure}\\
\vspace{0.2cm}
\begin{subfigure}[h]{0.3\textwidth}
  \centering
  \includegraphics[width=\linewidth]{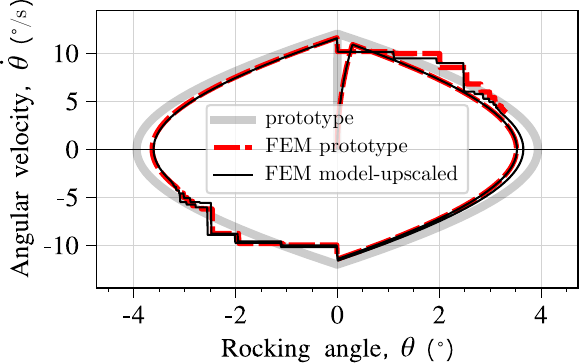}
	\caption{\footnotesize phase portrait, $W_{pa}=50$ kg}
\end{subfigure}\hspace{0.2cm}
\begin{subfigure}[h]{0.3\textwidth}
  \centering
  \includegraphics[width=\linewidth]{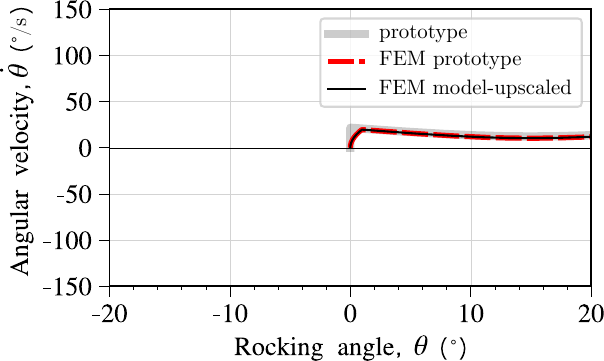}
	\caption{\footnotesize phase portrait, $W_{pb}=100$ kg}
\end{subfigure}
\caption{FE simulations of prototype and model responses ($\lambda=1/200$) and comparison with the solution of Eq. (\ref{eq:eom}). Both prototype and model are made of the same material.}
	\label{f:FEM_similar}
\end{figure*}

Figure \ref{f:FEM_similar} displays the FE results of the prototype and of the model, with the same material. The prototype and model responses coincide. It is worth noticing that, for both systems, the FE results agree remarkably well with the numerical solution of Eq. (\ref{eq:eom}). Indeed, whilst rocking, sliding and up-lifting happen together, the former prevails for slender blocks. Nevertheless, a shift in the rocking angle between the FE results and the rocking equation of motion appears after the first impact with the ground (at 1.3 ms from the arrival time). This is mainly due to repeated impacts with the ground that alter the motion when the rocking angle changes sign.

\subsubsection{Non-slender blocks}
\begin{figure*}[ht]
\centering
\begin{subfigure}[h]{0.3\textwidth}
  \centering
  \includegraphics[width=\linewidth]{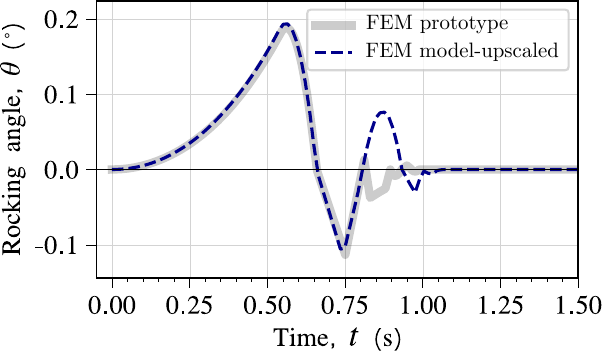}
	\caption{\footnotesize rocking angle.}
\end{subfigure}\hspace{0.2cm}
\begin{subfigure}[h]{0.3\textwidth}
  \centering
  \includegraphics[width=\linewidth]{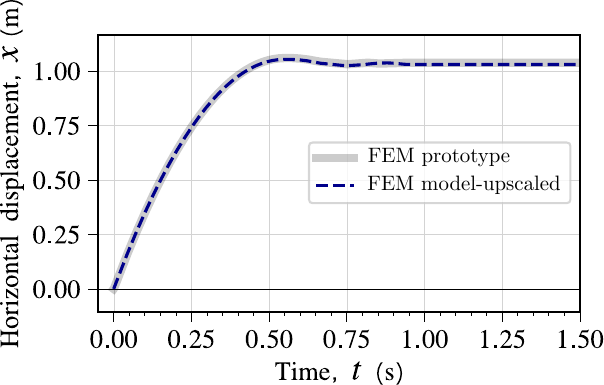}
	\caption{\footnotesize horizontal displacement.}
\end{subfigure} \\
\vspace{0.2cm}
\begin{subfigure}[h]{0.3\textwidth}
  \centering
  \includegraphics[width=\linewidth]{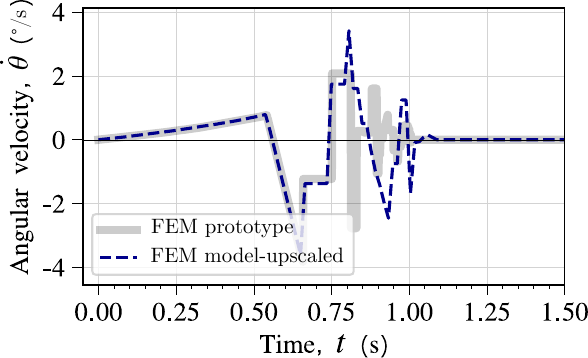}
	\caption{\footnotesize angular velocity}
\end{subfigure}\hspace{0.2cm}
\begin{subfigure}[h]{0.3\textwidth}
  \centering
  \includegraphics[width=\linewidth]{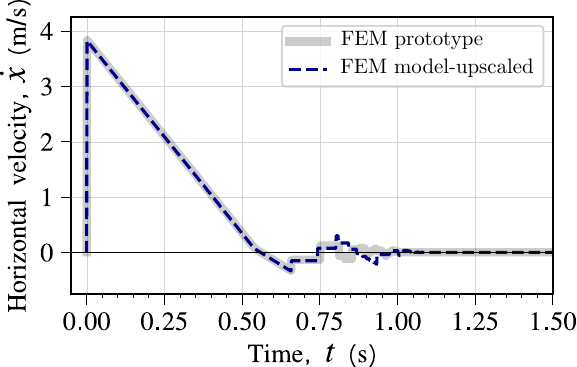}
	\caption{\footnotesize horizontal velocity.}
\end{subfigure}
\caption{FE simulations of prototype and model responses ($\lambda=1/200$) for non-slender blocks ($\alpha=35^{\circ}$) subjected to a prototype explosive charge of 2500 kg at 3 m. The prototype and model are made of the same material.}
	\label{f:FEM_nonslender}
\end{figure*}

\begin{figure*}[ht]
  \centering
  \includegraphics[width=\linewidth]{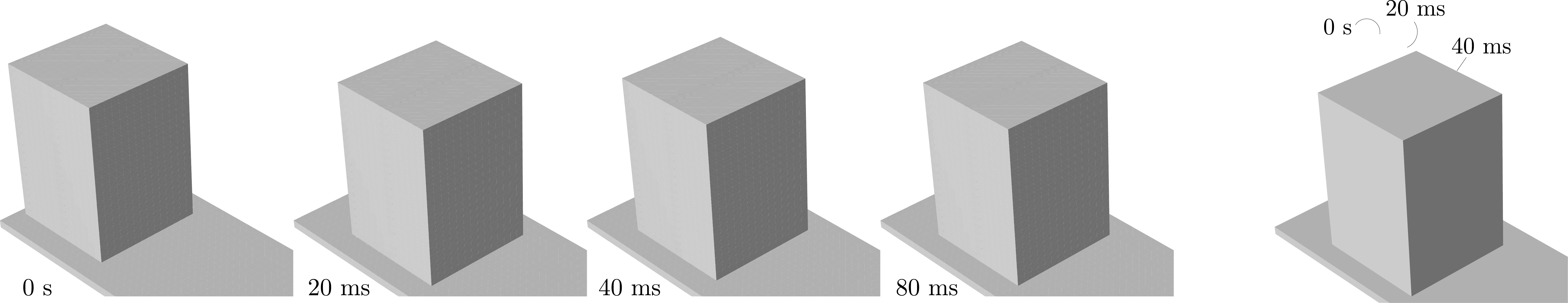}
\caption{Dynamic evolution of a non-slender rigid block subjected to $2500$ kg, at a stand-off distance of 3 m, according to the FE simulations (sliding is predominant).}
	\label{f:block_FE_sliding}
\end{figure*}

\noindent Sliding predominant responses are now investigated. We assume a non-slender prototype block, with height $2h_p=10$ m, slenderness angle $\alpha=35^{\circ}$, and depth and width $2w_p = 2b_p = 3.5$ m. The prototype model is subjected to an explosion of $2500$ kg of TNT, at $3$ m.\\
The model is characterized by a geometric scaling factor $\lambda = 1/200$ and unit density scaling factor. According to the scaling (Tab. \ref{tab:scaling_general}), we compute a similar explosive quantity equal to $14.26$ mg. The numerical results of the prototype and model systems are presented in Figure \ref{f:FEM_nonslender}, while the dynamic evolution of the prototype is displayed in Figure \ref{f:block_FE_sliding}. The system responses coincide perfectly up to $0.8$ s after the arrival of the shock wave. At this moment, the block is almost at rest (after sliding for a distance of approximately $1$ m), and a series of impacts with the ground base take place. Due to the successive impacts, the response of the model is found to differ from the one of the prototype. This might also be due to minor numerical issues related to contact calculations in the FE code. Nevertheless, the scaling laws assure similarity of the predominant quantity (sliding distance and velocity).

\clearpage
\section{Application to multi-drum masonry columns}
\label{sec:application}
\begin{figure*}[h]
\centering
\begin{subfigure}[h]{0.4\textwidth}
  \centering
  \includegraphics[width=0.8\linewidth]{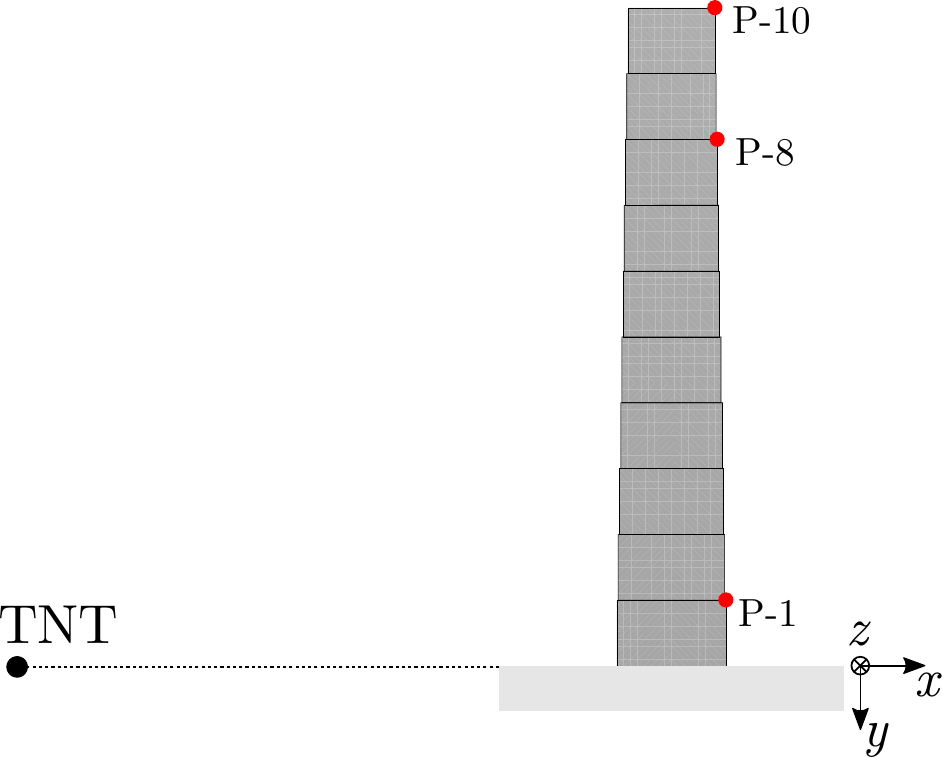}
	\caption{\footnotesize Multi-drum column with rectangular section.}
\end{subfigure}\hspace{0.2cm}
\begin{subfigure}[h]{0.4\textwidth}
  \centering
  \includegraphics[width=0.8\linewidth]{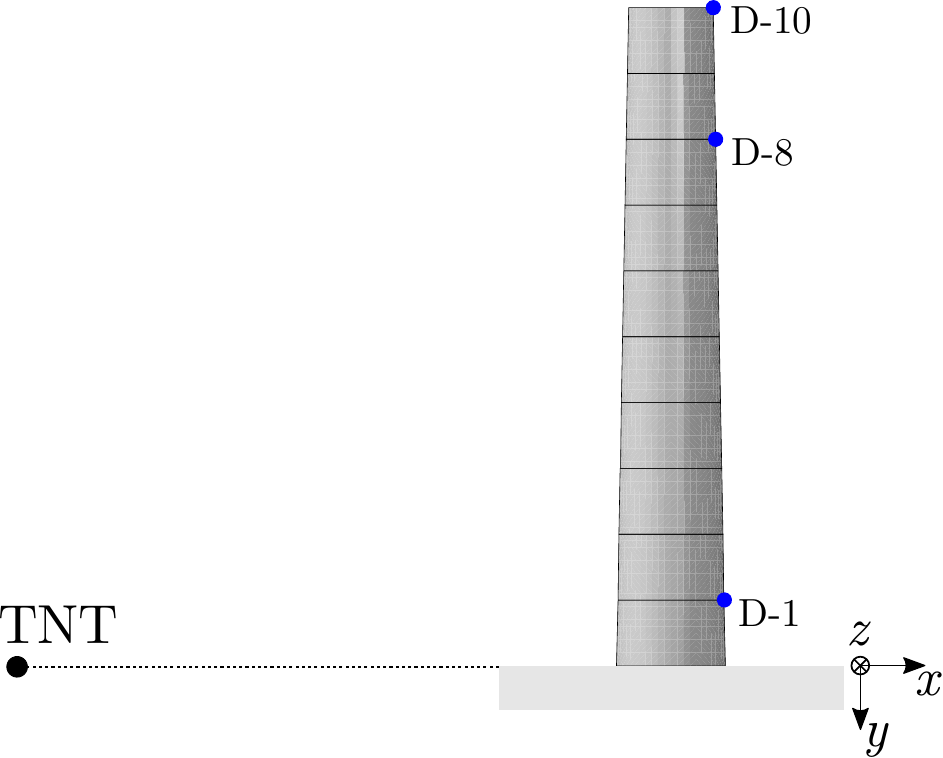}
	\caption{\footnotesize Multi-drum column with circular section.}
\end{subfigure} 
\caption{Geometry of the multi-drum columns with (a) square and (b) circular cross-section. Monitoring points are highlighted. Horizontal lines denote the contact of the drums.}
	\label{f:multidrum_model}
\end{figure*}

\noindent We consider multi-drum columns, made of fitted stone drums, placed on top of each other, avoiding the use of cement (mortar), see Figure \ref{f:multidrum_model}. These structures show rich dynamics \cite{psycharis,drumcolumns,pena,stefanou2011dynamic,makris2013planar,drosos2014shaking,sarhosis2016stability,fragiadakis2016vulnerability} and therefore they are a very interesting benchmark for the derived scaling laws.\\

Columns of two different sections are considered: one with a square section (Fig. \ref{f:multidrum_model}-a) and one of a circular one (Fig. \ref{f:multidrum_model}-b). The latter case renders the benchmark even more challenging due to wobbling and its dynamics, differing from rocking \cite{stefanou2011dynamic,STEFANOU20114325,vassiliou2018seismic,psycharis}. The multi-drum columns are $10$ m high and are composed of ten $1$ m high drums. The length of the square section and the diameter of the circular section vary from $1.65$, at the base, to $1.28$ m, at the top. The geometry is inspired by the external, Doric columns of the Parthenon \cite{neils2005parthenon}.\\

We stress that the scaling laws have been derived for masonry structures undergoing rigid-body motion. The following application shall show that the validity of the proposed scaling laws are also valid for multi-blocks structures.\\

Contrary to the previous validations with monolithic rigid blocks, we account for the deformability of the blocks. However, the scaling laws for the material response are not considered. Each drum is made of marble, which is considered to be (homogeneous) linearly elastic (with bulk and shear moduli equal to $50$ and $27$ GPa, respectively). For all numerical calculations, stresses are always found to be below the strength of the material. Deformations of the drums are, for all considered cases of blast actions, negligible. This shows that, for multi-block structures, as those studied here, a rigid-body response is predominant, as it is also the case for earthquake loads \cite{STEFANOU20114325,stefanou2015seismic,fragiadakis2016vulnerability}. We consider a Coulomb friction, with an angle of friction equal to $35^{\circ}$, at the interface of each drum. A hard contact formulation \cite{abaqus} is used.\\

\noindent
Blast loads are applied only on the columns' surfaces that are frontal to the explosive source considering two different methods:
\begin{itemize}
\item[(a)] \texttt{Method A}. Blast loads are applied relying on the best-fit interpolations shown in Appendix I, assuming simultaneity and uniformity of the pressure load.
\item[(b)] \texttt{Method B}. We account for the non-simultaneity of the load, the effects of surface rotation of the blocks, incident angle (Mach stem), and the relative distance between explosive and blocks. Following the approach in \citet{vannucci2017comparative,masi2020discrete} the position and the angle of incident of the center of the front surface of blocks are used to compute the blast loads.
\end{itemize}
In Appendix II, a third, more realistic blast load model (denoted with \texttt{Method C}) is considered. In this case, we account for all exposed surfaces (front, side, and rear ones) in the blast load calculation, which is the most realistic scenario.\\
We consider a model with same material as the prototype and a geometric scaling $\lambda=1/100$.
\subsection{Multi-drum column with square cross-section}
\begin{figure*}[h]
  \centering
  \begin{subfigure}[h]{1\textwidth}
  \centering
  \includegraphics[width=\linewidth]{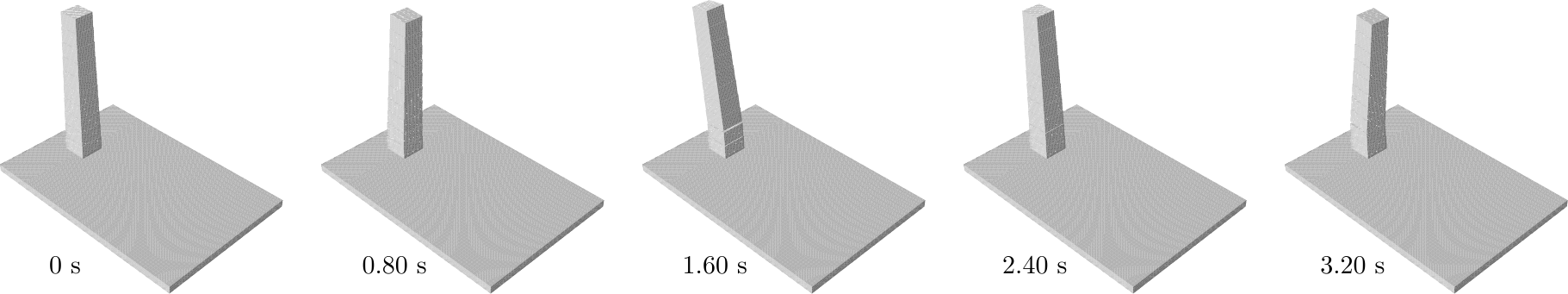}
	\caption{\footnotesize response mechanism under $250$ kg at $10$ m.}
\end{subfigure}\\
\vspace{0.2cm}
\begin{subfigure}[h]{1\textwidth}
  \centering
  \includegraphics[width=\linewidth]{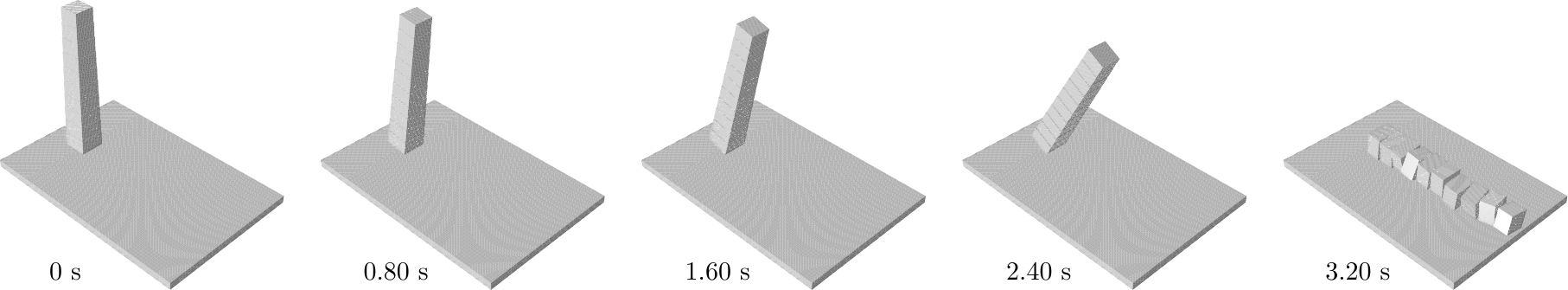}
	\caption{\footnotesize response mechanism under $500$ kg at $10$ m.}
\end{subfigure} 
\caption{Response mechanisms for a multi-drum column, with square cross-section, subjected to $250$ kg (a) and $500$ kg, at a stand-off distance of $10$ m. The time scale refers to the prototype system.}
	\label{f:multiprism_response}
\end{figure*}

\noindent The (prototype) columns are subjected to the loading arising from the denotation of (a) $250$ kg and (b) $500$ kg at a stand-off distance of $10$ m.
A schematic representation of the response mechanism to the two quantities of explosive is shown in Figure \ref{f:multiprism_response}.
\subsubsection{\texttt{Method A}}
\label{par:planar}
\begin{figure*}[h]
\centering
\begin{subfigure}[h]{0.3\textwidth}
  \centering
  \includegraphics[width=\linewidth]{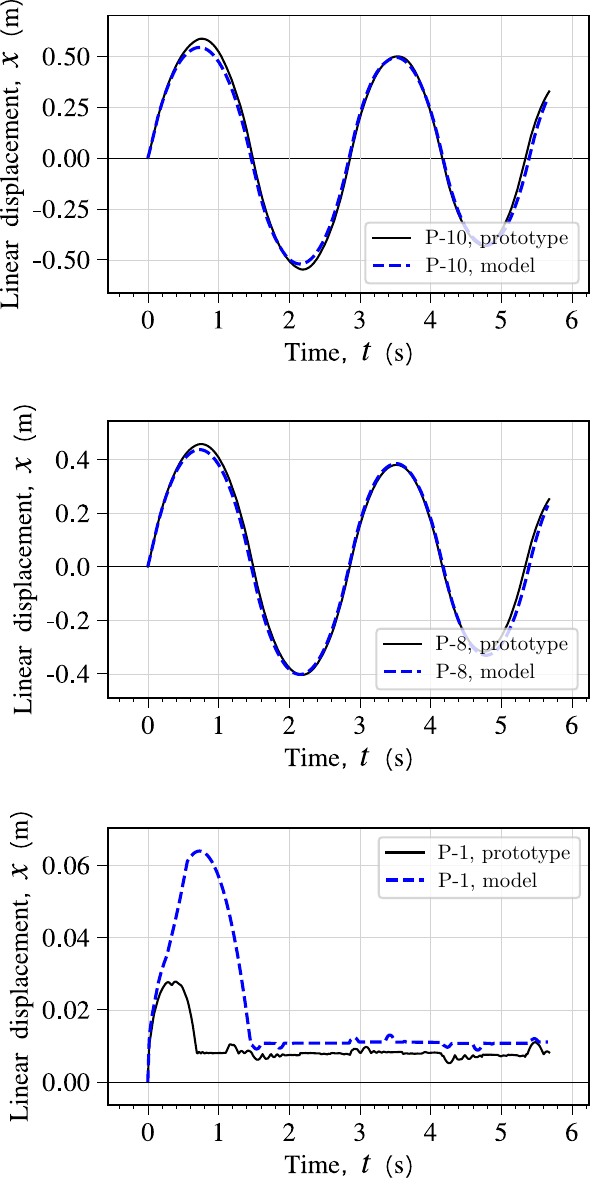}
	\caption{\footnotesize horizontal displacement.}
\end{subfigure}\hspace{0.5cm}
\begin{subfigure}[h]{0.3\textwidth}
  \centering
  \includegraphics[width=\linewidth]{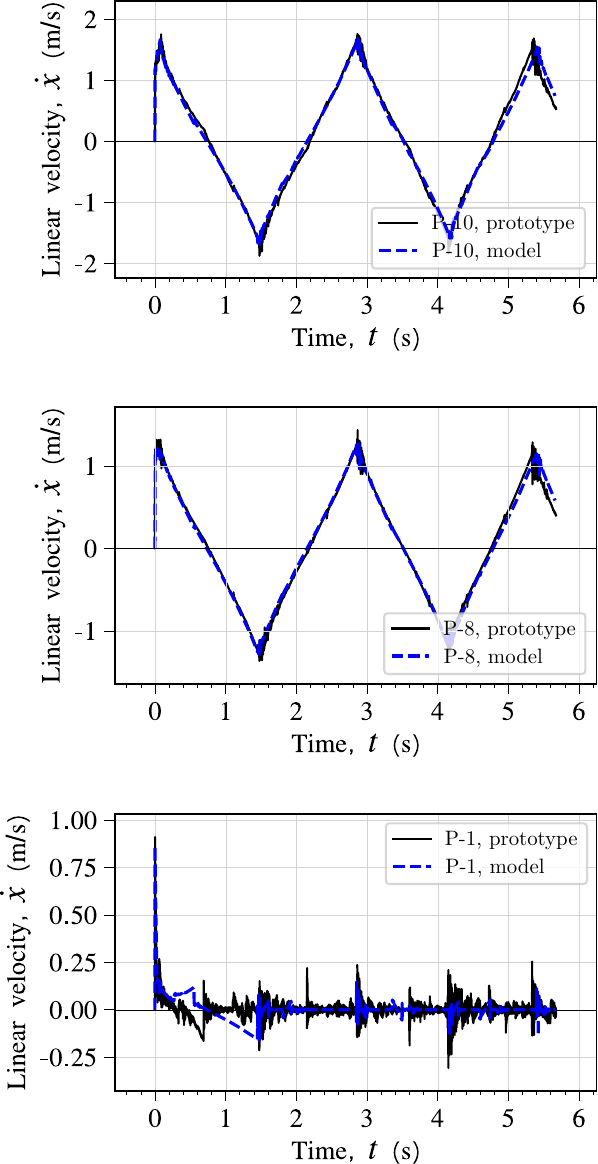}
	\caption{\footnotesize horizontal velocity.}
\end{subfigure} 
\caption{Comparison between the prototype response and the model response ($\lambda=1/100$), for a multi-drum column with square cross-section subjected to $250$ kg at $10$ m (\texttt{Method A}). Displacements and velocities of various monitoring points (cf. Fig. \ref{f:multidrum_model}) are represented in (a) and (b), respectively.}
	\label{f:prism_W250}
\end{figure*}

\begin{figure*}[h]
\centering
\begin{subfigure}[h]{0.3\textwidth}
  \centering
  \includegraphics[width=\linewidth]{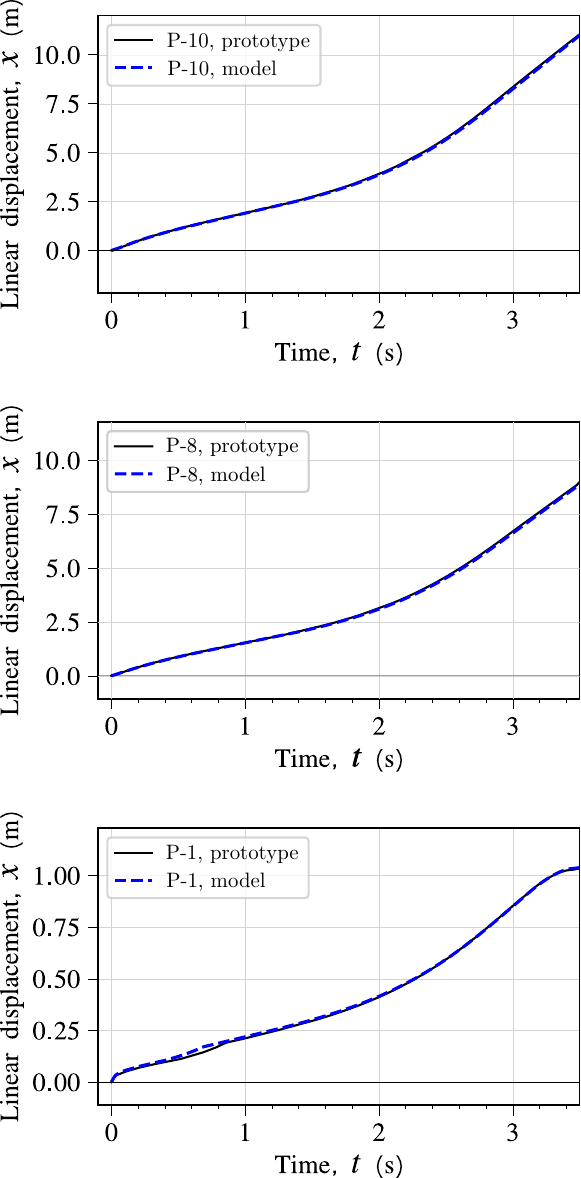}
	\caption{\footnotesize horizontal displacement.}
\end{subfigure}\hspace{0.5cm}
\begin{subfigure}[h]{0.3\textwidth}
  \centering
  \includegraphics[width=\linewidth]{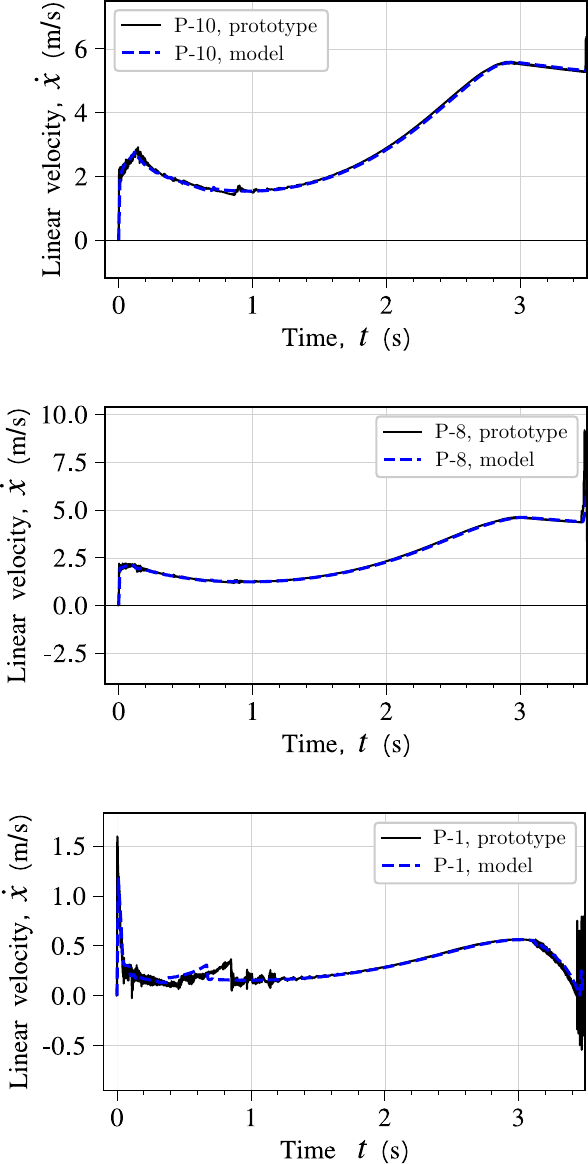}
	\caption{\footnotesize horizontal velocity.}
\end{subfigure} 
\caption{Comparison between the prototype response and the model response ($\lambda=1/100$), for a multi-drum column with cross-square section subjected to $500$ kg at $10$ m (\texttt{Method A}).}
	\label{f:prism_W500}
\end{figure*}

\noindent By relying on the simplified approach of considering a planar shock wave impinging only the front surface of the structure, we investigate the dynamic response and the validity of the proposed scaling laws for deformable multi-block masonry structures.\\
Figure \ref{f:prism_W250} displays the horizontal displacement and velocity (along $x$-axis, see Fig. \ref{f:multidrum_model}) at different monitoring points, for the column with square cross-section, subjected to an explosive quantity of $250$ kg at $10$ m away. The response of both the prototype and the model (up-scaled) are plotted. Figure \ref{f:prism_W500} refers instead to an explosive quantity of $500$ kg. In both cases, we can see that the overall response of the prototype is well captured by the model. This is particularly true for the response of the blocks at the top (i.e., P-5 to P-10, cf. Fig. \ref{f:multidrum_model}).\\

Differences, between the prototype and the model, are visible at the lower blocks (P-1 and P-3). This is mainly due to numerical errors and discrepancies related to the contact algorithm used. Furthermore, we only considered the scaling laws for the rigid-body motion. As a result, similarity for elastic deformation mechanisms is not guaranteed. The repeated impacts (exciting the algorithmic elastic response of the interfaces) are the main cause of the differences between the two systems. It is worth noticing that the same issue is found also for single rigid blocks. Nevertheless, the aforementioned comparisons proved that the elastic deformation at the interfaces has very small influence on the dynamic response. The deformations of the drums may, as well, have some minor influence on the differences between the model and the prototype.

\subsubsection{\texttt{Method B}}
\label{par:nonplanar}
\begin{figure*}[h]
\centering
\begin{subfigure}[h]{0.3\textwidth}
  \centering
  \includegraphics[width=\linewidth]{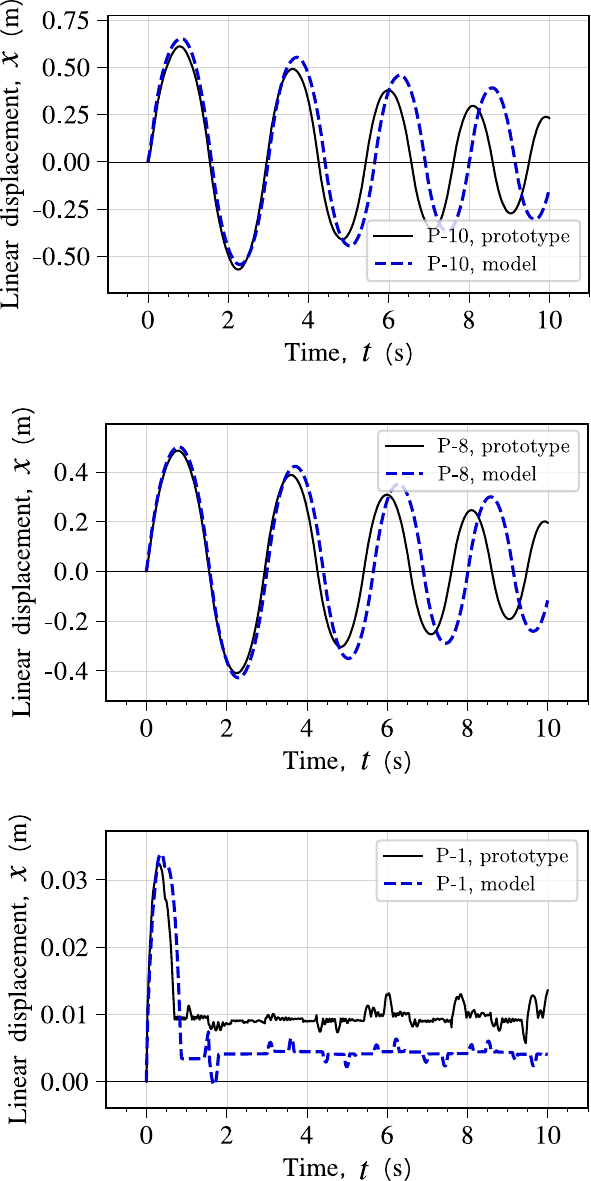}
	\caption{\footnotesize horizontal displacement.}
\end{subfigure}\hspace{0.2cm}
\begin{subfigure}[h]{0.3\textwidth}
  \centering
  \includegraphics[width=\linewidth]{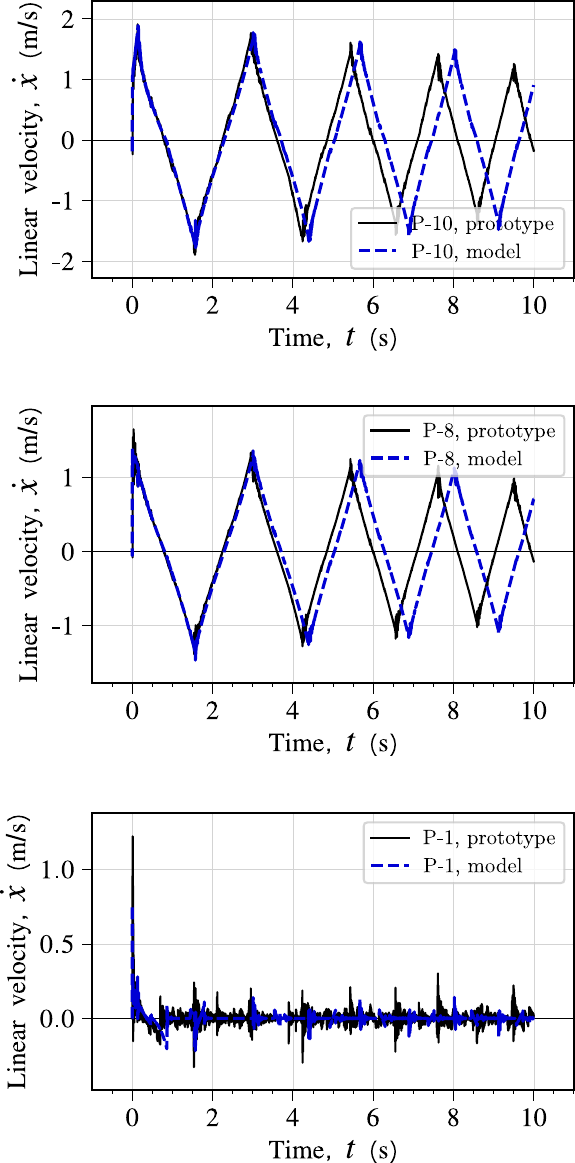}
	\caption{\footnotesize horizontal velocity.}
\end{subfigure} 
\caption{Comparison between the prototype response and the model response ($\lambda=1/100$), for a multi-drum column with cross-square section subjected to $250$ kg at $10$ m (\texttt{Method B}).}
	\label{f:FEM_prism_WA_B}
\end{figure*}

\begin{figure*}[h]
\centering
\begin{subfigure}[h]{0.3\textwidth}
  \centering
  \includegraphics[width=\linewidth]{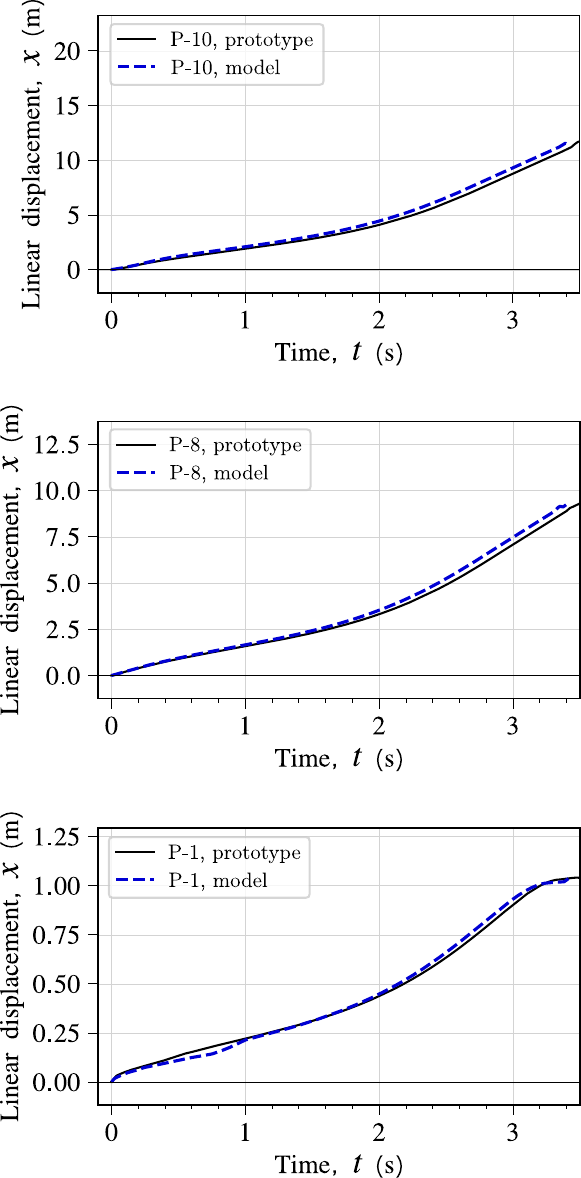}
	\caption{\footnotesize horizontal displacement.}
\end{subfigure}\hspace{0.2cm}
\begin{subfigure}[h]{0.3\textwidth}
  \centering
  \includegraphics[width=\linewidth]{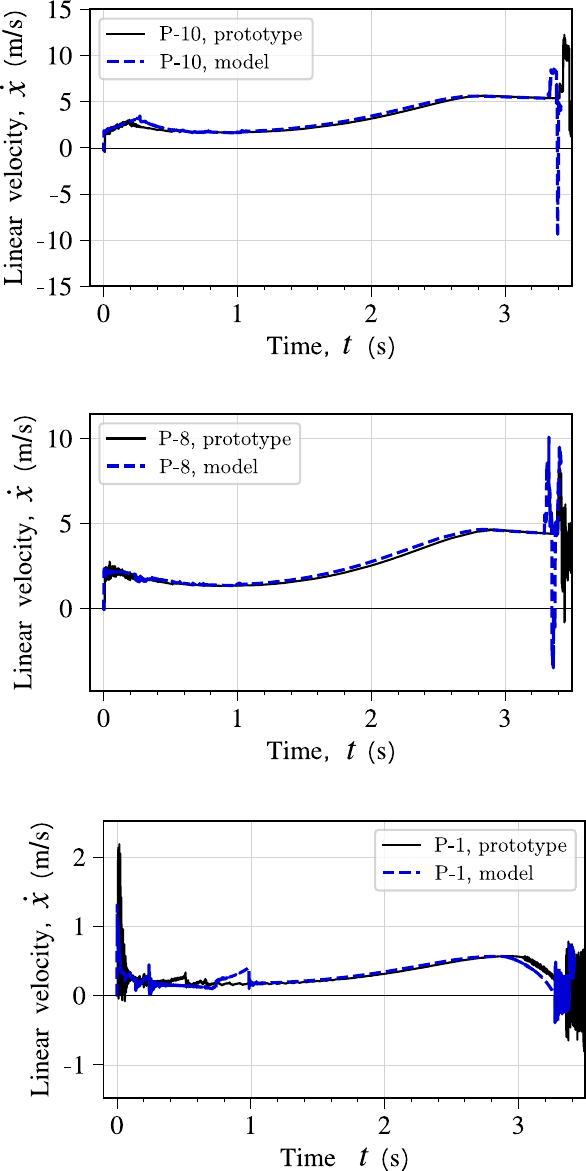}
	\caption{\footnotesize horizontal velocity.}
\end{subfigure} 
\caption{Comparison between the prototype response and the model response ($\lambda=1/100$), for a multi-drum column with cross-square section subjected to $500$ kg at $10$ m (\texttt{Method B}).}
	\label{f:FEM_prism_WB_B}
\end{figure*}

\noindent The more detailed characterization of blast loads is here adopted and the scaling laws are tested. In particular, we consider the spatial and temporal effects of an hemispherical shock wave. More complex phenomena than those considered in deriving the similarity laws are thus accounted for. In particular, the shock wave velocity (or, equivalently, the fact that the blast wave does not impinge all the front surface simultaneously) represents a main factor for testing the performance of the derived scaling laws.\\
Figures \ref{f:FEM_prism_WA_B} and \ref{f:FEM_prism_WB_B} show the response of the prototype and (upscaled) model against $250$ and $500$ kg prototype explosive charges, respectively. It is worth noticing that the main response of the system is well predicted by the model. Nevertheless, some differences exist. This is particularly true for the bottom blocks. Besides numerical issues related to the contact algorithm of the FE code used, the discrepancy can be also owed to the fact that the shock wave velocity, i.e., the shifting of the arrival time for each surface, is not taken into account by the proposed scaling laws. This results to some differences between the model and the prototype, but they are minor (see Fig.s \ref{f:FEM_prism_WA_B} and \ref{f:FEM_prism_WB_B}).

\subsection{Multi-drum column with circular cross-section}
\begin{figure*}[h]
  \centering
\begin{subfigure}[h]{1\textwidth}
  \centering
  \includegraphics[width=\linewidth]{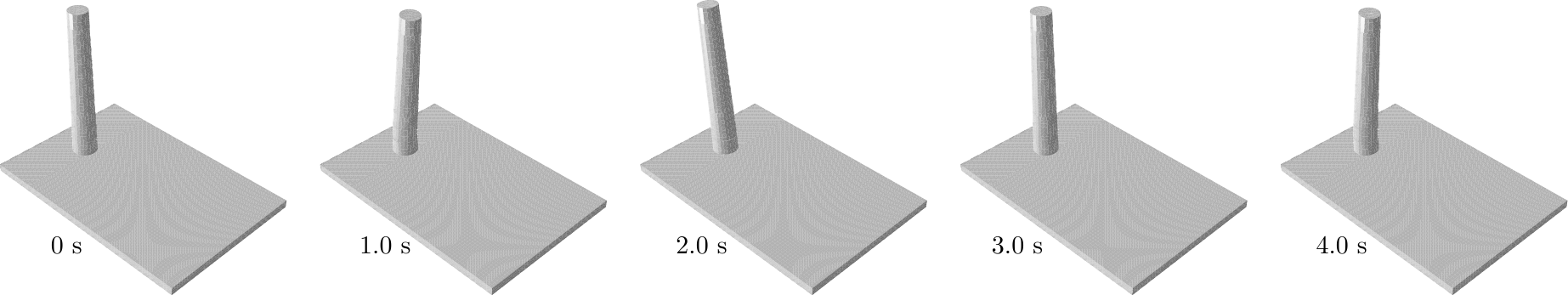}
	\caption{\footnotesize response mechanism under $200$ kg at $10$ m.}
\end{subfigure}\\
\vspace{0.2cm}
\begin{subfigure}[h]{1\textwidth}
  \centering
  \includegraphics[width=\linewidth]{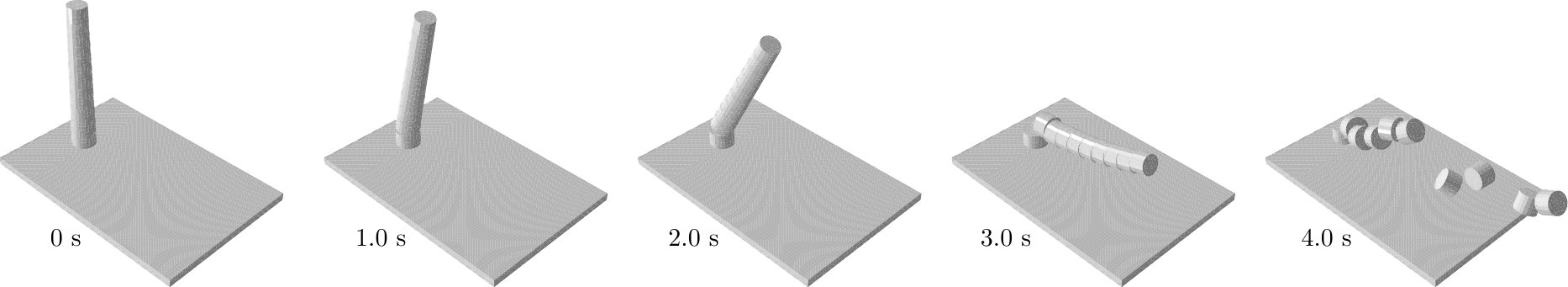}
	\caption{\footnotesize response mechanism under $400$ kg at $10$ m.}
\end{subfigure} 
\caption{Response mechanisms for a multi-drum column, with circular cross-section, subjected to $200$ kg (a) and $400$ kg, at a stand-off distance of $10$ m. Rocking and wobbling co-exist.}
	\label{f:multicirc_response}
\end{figure*}

\noindent The (prototype) columns are subjected to the loading arising from the detonation of (a) $200$ kg and (b) $400$ kg at a stand-off distance of $10$ m.
A schematic representation of the response mechanism under the two quantities of explosive is shown in Figure \ref{f:multicirc_response}.

\subsubsection{\texttt{Method A}}
\label{par:planarS}
\begin{figure*}
\centering
\begin{subfigure}[h]{0.3\textwidth}
  \centering
  \includegraphics[width=\linewidth]{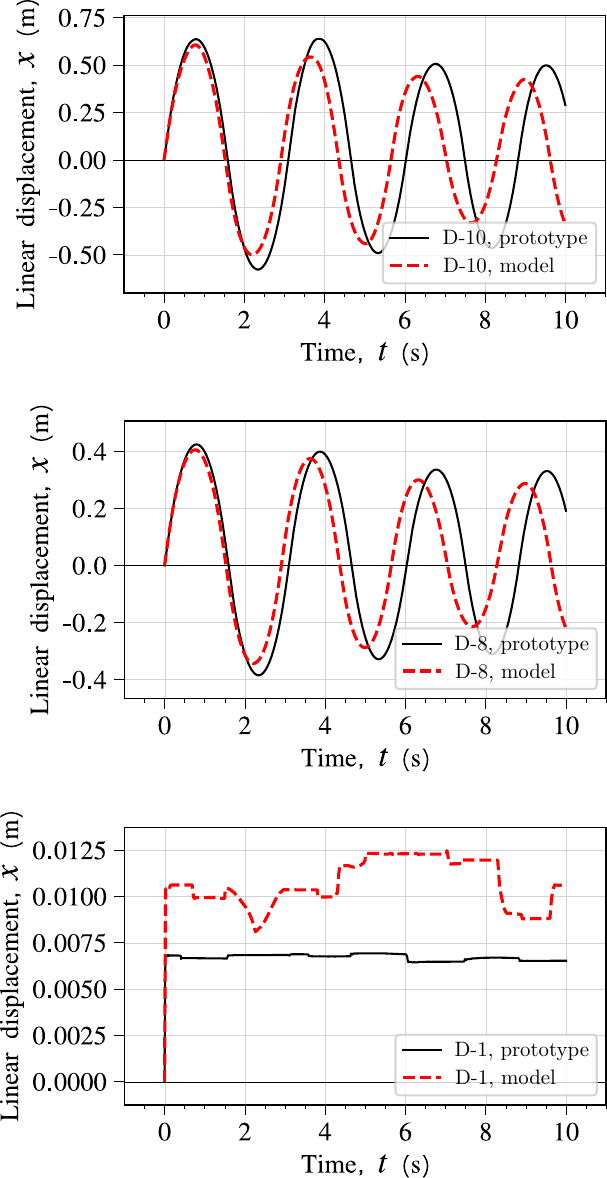}
	\caption{\footnotesize horizontal displacement.}
\end{subfigure}\hspace{0.2cm}
\begin{subfigure}[h]{0.3\textwidth}
  \centering
  \includegraphics[width=\linewidth]{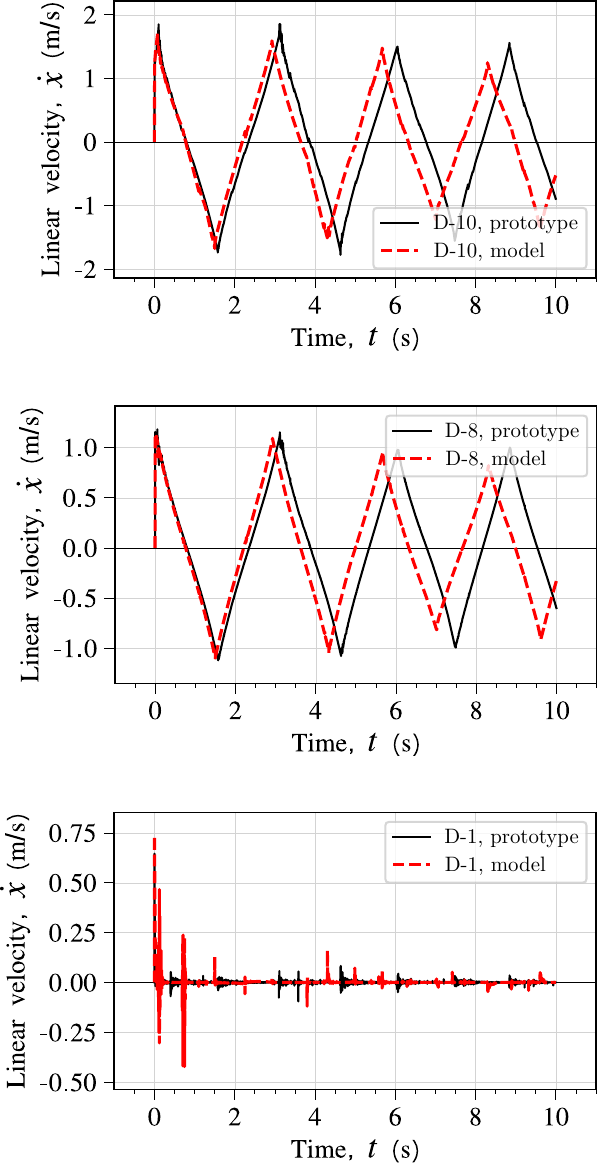}
	\caption{\footnotesize horizontal velocity.}
\end{subfigure} 
\caption{Comparison between the prototype response and the model response ($\lambda=1/100$), for a multi-drum column with circular cross-section subjected to $200$ kg at $10$ m (\texttt{Method A}). Displacements and velocities of various monitoring points (cf. Fig. \ref{f:multidrum_model}) are represented in (a) and (b), respectively.}
	\label{f:circ_W250}
\end{figure*}

\begin{figure*}
\centering
\begin{subfigure}[h]{0.3\textwidth}
  \centering
  \includegraphics[width=\linewidth]{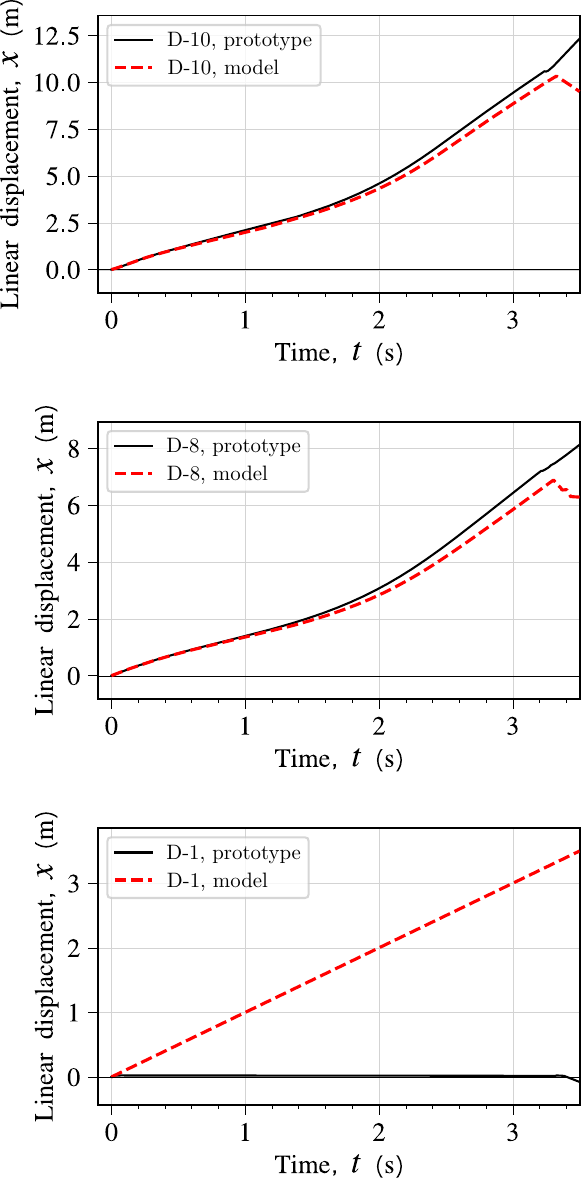}
	\caption{\footnotesize horizontal displacement.}
\end{subfigure}\hspace{0.2cm}
\begin{subfigure}[h]{0.3\textwidth}
  \centering
  \includegraphics[width=\linewidth]{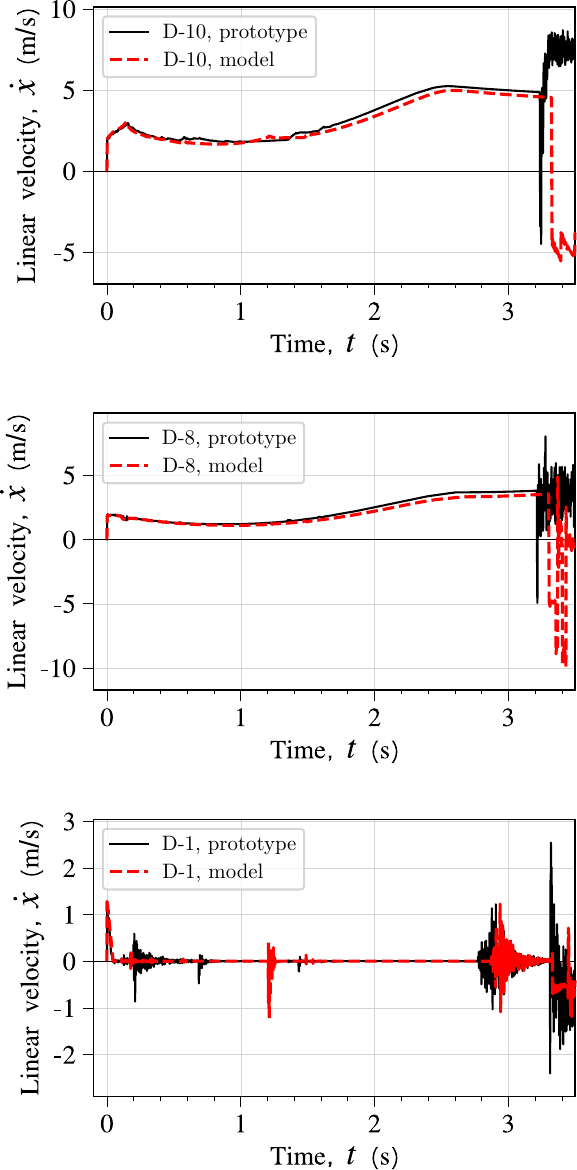}
	\caption{\footnotesize horizontal velocity.}
\end{subfigure} 
\caption{Comparison between the prototype response and the model response ($\lambda=1/100$), for a multi-drum column with circular cross-section subjected to $400$ kg at $10$ m (\texttt{Method A}).}
	\label{f:circ_W500}
\end{figure*}

\noindent Assuming a planar shock wave impinging the front surfaces at the same time, we study the response of a multi-drum column with a circular cross-section. Figure \ref{f:circ_W250} and \ref{f:circ_W500} display the displacement and velocity, at several monitoring points, of a column with circular cross-section under $200$ kg and $400$ kg of TNT equivalent, respectively. Once more, the first peak response is extremely well captured by the model. Nevertheless, after the first impact ($\approx 1.5$ s), the response predicted by the model is found to diverge from that of the prototype. The reason relies on the complex response of the system, characterized by repeated impacts, wobbling motion, stick-slips, and rocking, see e.g. \citet{stefanou2011dynamic,vassiliou2018seismic,STEFANOU20114325,psycharis}. Examining carefully the results, we observe that the model does not capture exactly the period of the movement. By considering the full equation of motion for wobbling \cite{stefanou2011dynamic} one can derive the required, exact scaling laws for columns of circular cross-section. However, this exceeds the scope of this work. For instance, notice that this is not the case for a square column, for which wobbling does not take place.\\
Nevertheless, even for the case of circular columns, the scaling laws give satisfactory results and accurately predict the first peak response, as well as the main response mechanism (collapse against $400$ kg of TNT).\\
In Figure \ref{f:multicirc_response} we present the evolution of the response of the structure due to the two quantities of explosive.

\subsubsection{\texttt{Method B}}
\label{par:nonplanarS}
\begin{figure*}[h]
\centering
\begin{subfigure}[h]{0.3\textwidth}
  \centering
  \includegraphics[width=\linewidth]{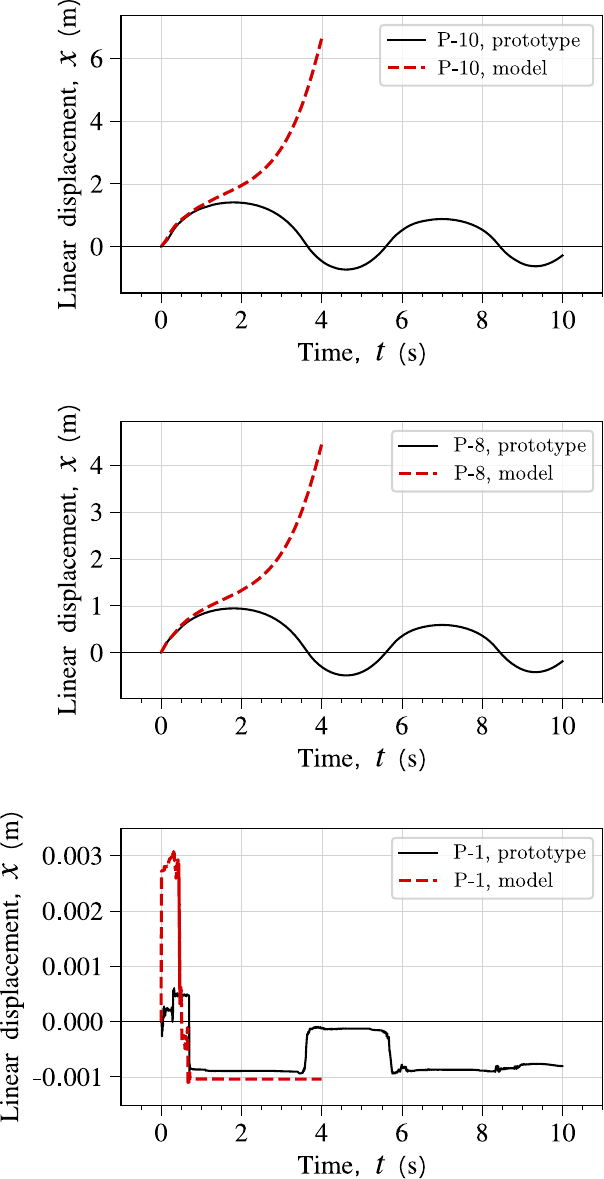}
	\caption{\footnotesize horizontal displacement.}
\end{subfigure}\hspace{0.2cm}
\begin{subfigure}[h]{0.3\textwidth}
  \centering
  \includegraphics[width=\linewidth]{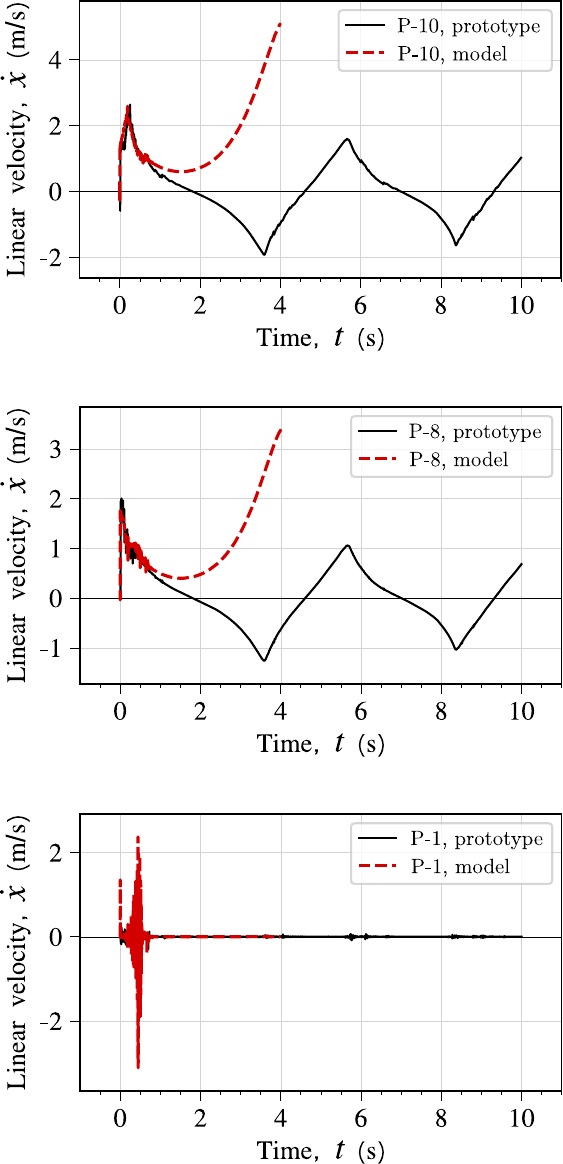}
	\caption{\footnotesize horizontal velocity.}
\end{subfigure} 
\caption{Comparison between the prototype response and the model response ($\lambda=1/100$), for a multi-drum column with circular cross-section subjected to $200$ kg at $10$ m (\texttt{Method B}).}
	\label{f:FEM_circ_WA_B}
\end{figure*}

\noindent We consider here the response of a circular cross-section column by using the more detailed description of the blast load. Contrary to the case of the column with square section, the responses of the prototype and of the model differ significantly for the load scenario of $200$ kg. Figure \ref{f:FEM_circ_WA_B} displays the horizontal displacement and velocities at several monitoring points. The very first instants after the shock arrival ($\approx 1$ s) are well captured by the model. Nevertheless, the model overturns while the prototype does not.\\

\begin{figure*}[h]
\centering
\begin{subfigure}[h]{0.3\textwidth}
  \centering
  \includegraphics[width=\linewidth]{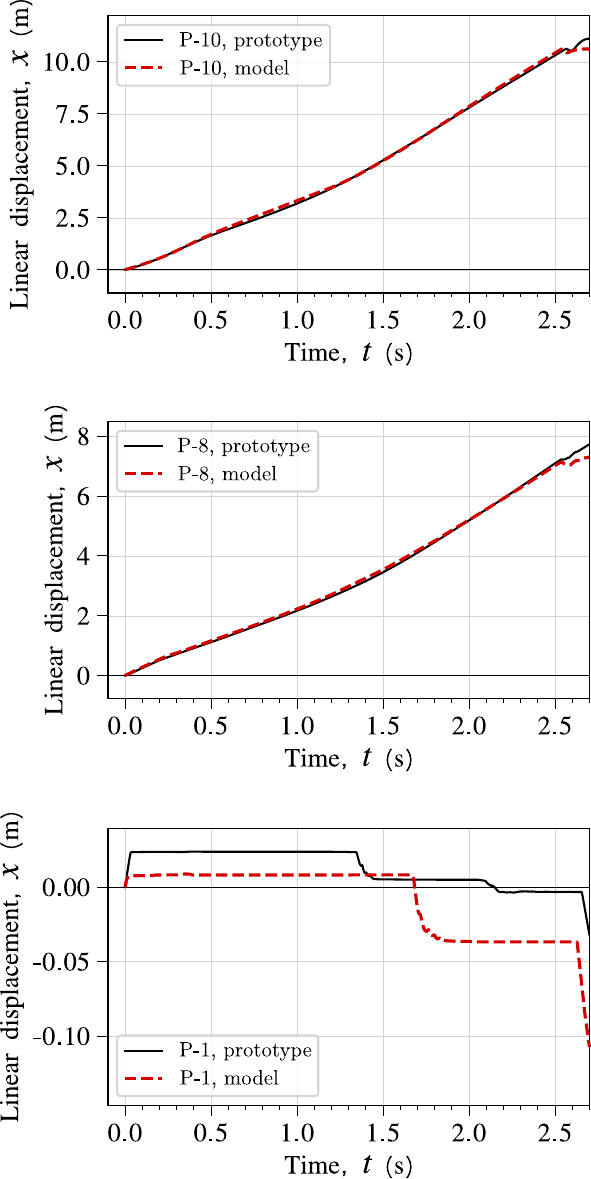}
	\caption{\footnotesize horizontal displacement.}
\end{subfigure}\hspace{0.2cm}
\begin{subfigure}[h]{0.3\textwidth}
  \centering
  \includegraphics[width=\linewidth]{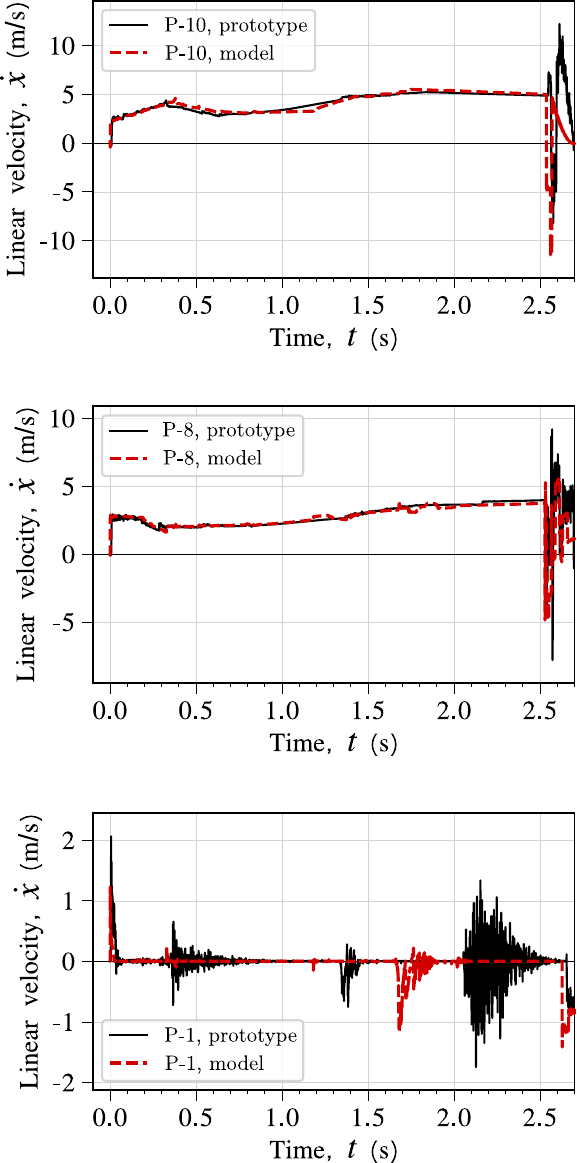}
	\caption{\footnotesize horizontal velocity.}
\end{subfigure} 
\caption{Comparison between the prototype response and the model response ($\lambda=1/100$), for a multi-drum column with circular cross-section subjected to $400$ kg at $10$ m (\texttt{Method B}).}
	\label{f:FEM_circ_WB_B}
\end{figure*}

The difference relies on the complex  dynamics of wobbling which, enhanced by the effect of dissimilar shock wave velocity (between prototype and model), causes important differences. Further investigations and developments to include the effects of differing shock wave arrival time are needed. Nevertheless, for the detonation of $400$ kg, the model, as expected, perfectly predicts the prototype response as overturning takes place (see Fig. \ref{f:FEM_circ_WB_B}). Indeed, under such a high quantity of explosive, rocking is dominant and the shock wave velocity is far higher \cite{vannucci2017comparative}. Consequently the shifting of the arrival time of the shock front on the impinged surfaces has minor effects.

\section{Concluding remarks}
\noindent Experiments of masonry structures under blast loads are rare in the available literature.\\
This is an obstacle for designing appropriate protective devices to preserve masonry structures against blasts.
Indeed, experimental testing of this kind of elements is particularly challenging due to the complex structural dynamic response of masonry. Furthermore, field testing shows several limitations related to cost, environmental hazards, safety risks, and repeatability. In addition, one of main restrictions to full-scale field testing stems from the explosive charge weight, which cannot exceed safety values (for the personnel and the monitoring equipment). Field testing of reduced scales prototype is thus necessary. In order to perform and accurately design scaled tests, scaling laws for both the blast loads and the specimens are mandatory.

We aim at providing scaling laws for masonry structures subjected to explosions. Based on previous works \cite{rockingmasi,masi2020resistance}, we derived similarity laws for the rigid-body motion of monolithic and blocky-masonry structures, considering empirical models for the blast actions. In contrast with the well-known Hopkinson-Cranz scaling laws \cite{hopki,cranz}, the proposed scaling laws allow to design experiments by reducing the blast intensity, which is compelling for safe experiments.

The proposed scaling laws were validated against numerical cases of monolithic prototypes and models and through three-dimensional Finite Element simulations. 
Finally, multi-drum columns, typical examples of key load-carrying elements in ancient masonry structures and monuments, were investigated. In particular, we showed that the scaling laws are valid for multi-block, deformable structures.\\

We give first insights of how reduced-scale experiments of ancient and modern masonry structures can be designed. Further investigations including richer dynamics, such as the wobbling motion, and the consideration of the blast wave front propagation along the structure are needed. Nevertheless, the derived scaling laws can be directly used for the design of preliminary experiments of masonry structures.

\subsection{Acknowledgments}
The authors would like to acknowledge Professor Paolo Vannucci (LMV, UMR 8100, Universit\'{e} de Versailles et Saint-Quentin) for various discussions relating to dynamics and blast loading of structures.

\pagebreak
\appendix
\section{Appendix A}
\begin{figure}[h]
	\centering
	\includegraphics[width=0.75\textwidth]{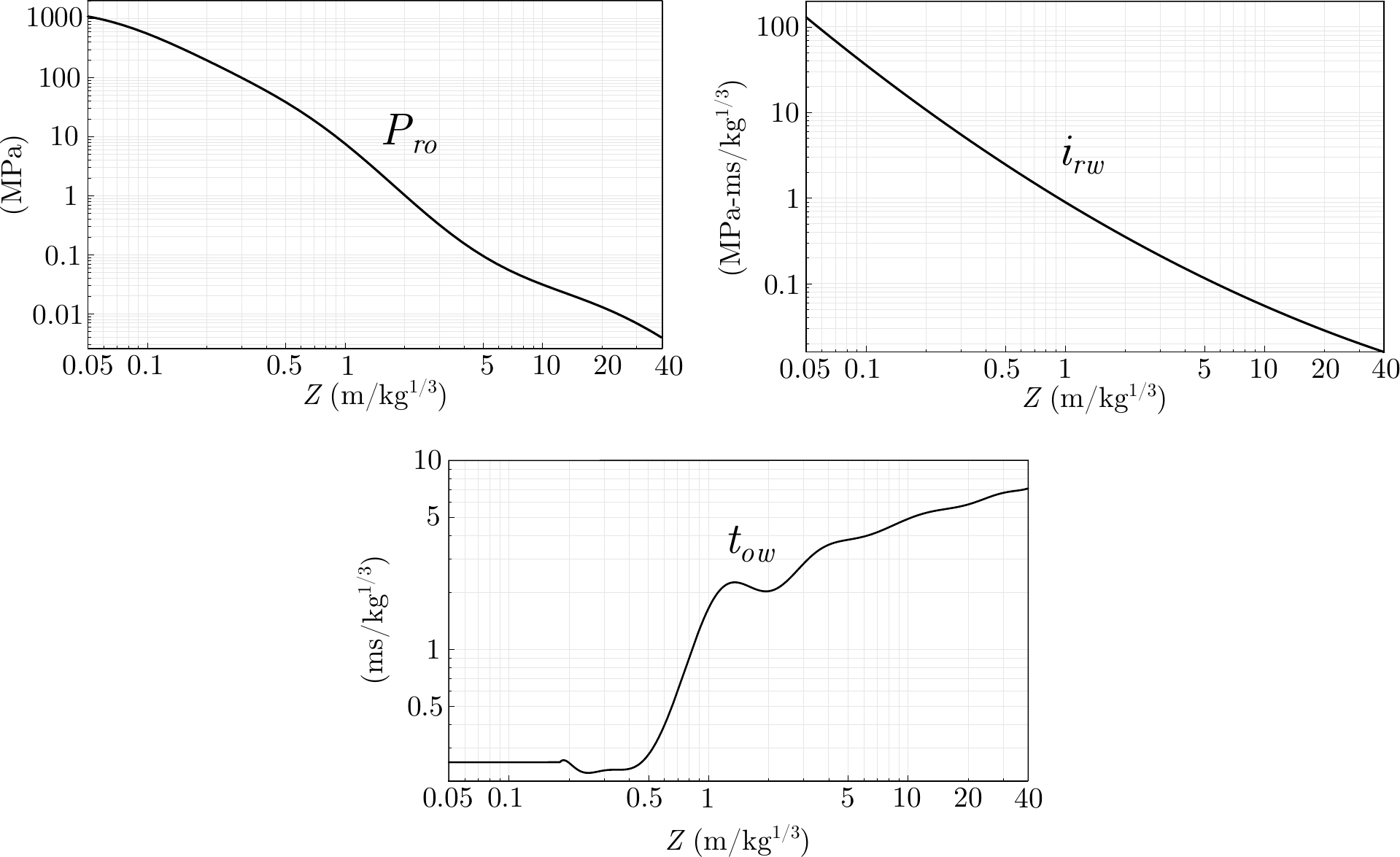}
	\caption{Analytical interpolations for blast loading as functions of the scaled distance, $Z$: reflected pressure, scaled reflected impulse and positive duration.}
	\label{f:appinterpolations}
\end{figure}

The expressions for the air-blast parameters for a surface burst are given below and presented in Figure \ref{f:appinterpolations}. For more details we refer to \citet{vannucci2017comparative,rockingmasi}.\\
\begin{itemize}
	\item normal reflected pressure peak $P_{ro}$:
		\begin{eqnarray*}
		P_{ro}(Z) & = & \Big(1+\frac{1}{2 e^{10Z}}\Big) \exp\Big[2.0304-1.8036\ln Z\\
		&&-0.09293\ln^2Z -0.8779\sin(\ln Z)-0.3603\sin^2(\ln Z)\Big]
		\end{eqnarray*}
	\item scaled and effective positive reflected impulse $i_{rw}$, $i_{r}$:
		\begin{eqnarray*}
		i_{rw}(Z,W)&=&\exp\Big[-0.110157-1.40609\ln Z+0.0847358\ln^2Z\Big],\\
		i_{r}(Z)&=& W^{\frac{1}{3}}i_{rw}(Z,W)
		\end{eqnarray*}
	\item  scaled and effective arrival time $t_{Aw}$, $t_{A}$:
		\begin{eqnarray*}
		t_{Aw}(Z,W)& =&\begin{cases}&
			\exp\Big[-0.6847+1.4288\ln Z+0.0290\ln^2Z \\
			&\qquad \hspace{0.2cm}+0.4108\sin(\ln Z)\Big] \qquad \text{if } Z \geq 0.18 \text{ m/kg}^{1/3},\\
		&0.0315495 \hspace{3cm} \text{if } Z < 0.18 \text{ m/kg}^{1/3},\\
		\end{cases}\\ \vspace{0.5cm}\\
		t_{A}(Z)&=&W^{\frac{1}{3}}t_{Aw}(Z,W)
		\end{eqnarray*}
	\item scaled and effective positive duration time $t_{ow}$, $t_{o}$:
		\begin{eqnarray*}
		t_{ow}(Z,W)&=&\begin{cases}&\begin{split}
			  \exp \bigg[&0.592+2.913 \ln Z-1.287 \ln^2 Z-1.788\ln^3Z\\
			  & +1.151\ln^4Z+0.325\ln^5Z -0.383\ln^6Z\\ &+0.090\ln^7Z-0.004\ln^8Z-0.0004\ln^9Z\\ &+0.537\cos^7\Big[1.032\left(\ln Z-0.859\right)\Big]\\
			  &\quad  \sinh\Big[1.088\left(\ln Z-2.023\right)\Big] \bigg]
			\end{split}\\
		&\hspace{5cm}\text{if } Z \geq 0.18 \text{ m/kg}^{1/3},\\
		& 0.251703  \hspace{3.7cm}  \text{if } Z < 0.18 \text{ m/kg}^{1/3},\\
		\end{cases}\\ \vspace{0.3cm}\\
		t_{o}(Z) &= & W^{\frac{1}{3}}t_{ow}(Z,W)
		\end{eqnarray*}

\end{itemize}

\section{Appendix B}
The scaling laws for a (deformable) multi-drum column, with square cross-section, are investigated by considering a more realistic characterization of the blast loading, denoted here with \texttt{Method C}. In particular, we consider here all the exposed surfaces of the structure (front, top, lateral, and rear), and not only the surface frontal to the explosive source (cf. Figure \ref{fig:scheme}). For each surface, we account for the non-simultaneity of the load, the effects of surface rotation of the blocks, incident angle (Mach stem), and the relative distance between explosive and blocks, as for \texttt{Method B}. Following the approach in \cite{vannucci2017comparative,masi2020discrete} the position and the angle of incident of the front surface of blocks are used to computed the blast loads.\\

The prototype is subjected to several explosive quantities: $250$, $500$, $750$, and $1000$ kg. Differently from the above results, the overturning happens for a quantity of $1000$ kg. Indeed, the effects of the pressure acting on rear surfaces are stabilizing \cite{rockingmasi}.\\

\begin{figure*}[htb]
\centering
\begin{subfigure}[b]{0.63\textwidth}
  \centering
  \includegraphics[width=\linewidth]{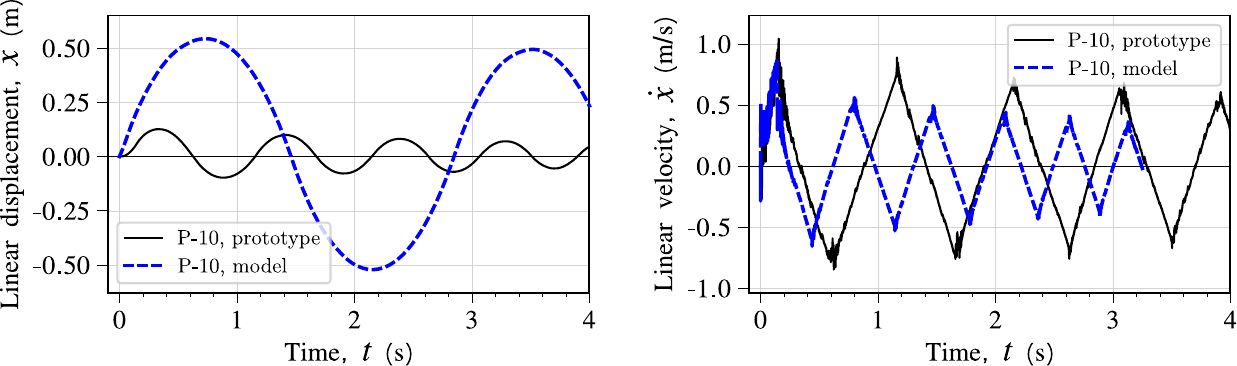}
	\caption{\footnotesize 250 kg}
	\label{f:meth_C1}
\end{subfigure}\\
\begin{subfigure}[b]{0.63\textwidth}
  \centering
  \includegraphics[width=\linewidth]{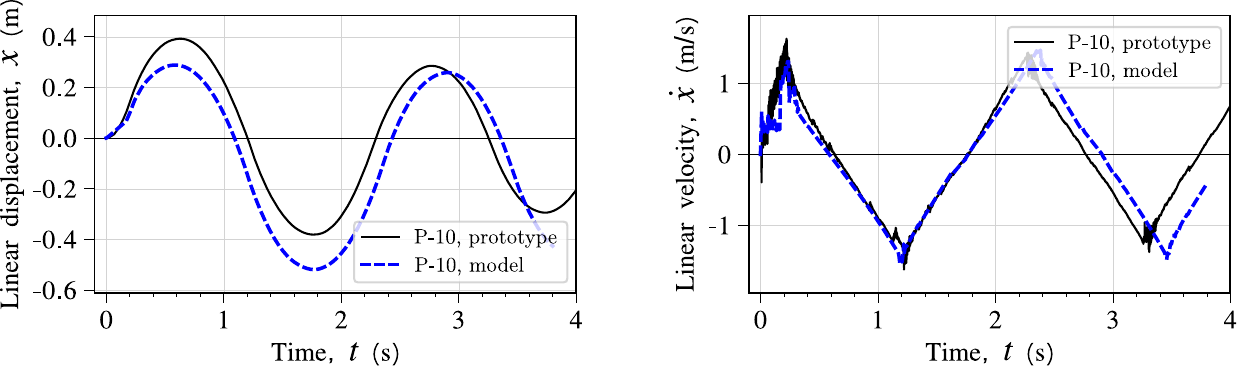}
	\caption{\footnotesize 500 kg}
	\label{f:meth_C2}
\end{subfigure}\\
\begin{subfigure}[b]{0.63\textwidth}
  \centering
  \includegraphics[width=\linewidth]{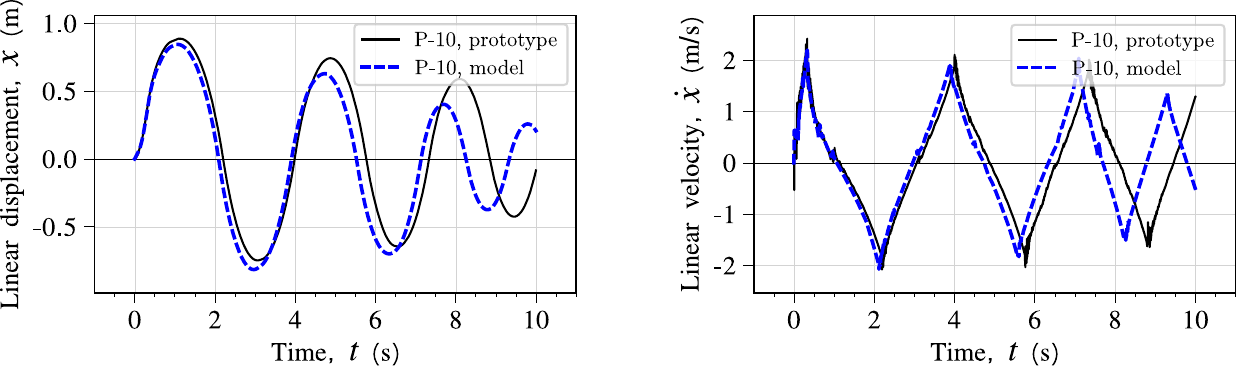}
	\caption{\footnotesize 750 kg}
	\label{f:meth_C3}
\end{subfigure}\\
\begin{subfigure}[b]{0.63\textwidth}
  \centering
  \includegraphics[width=\linewidth]{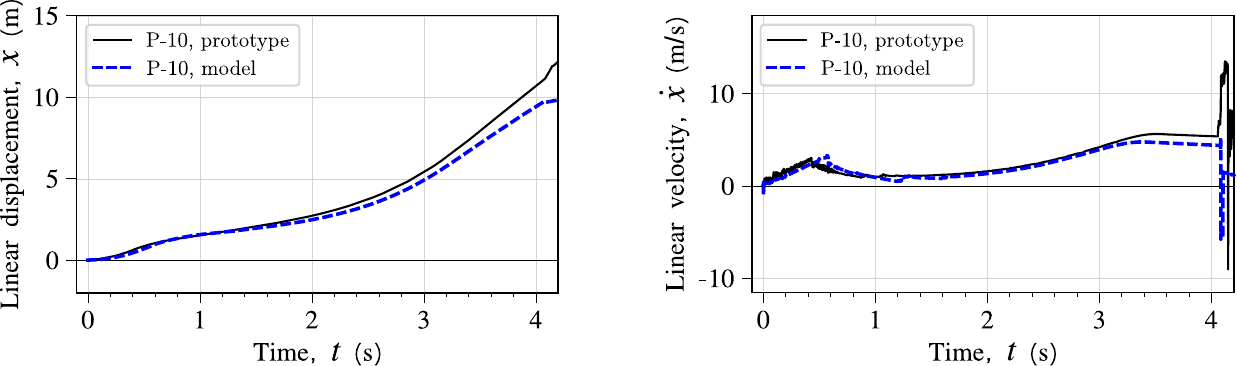}
	\caption{\footnotesize 1000 kg}
	\label{f:meth_C4}
\end{subfigure}
\caption{Comparison between the prototype response and the model response ($\lambda=1/100$), for a multi-drum column with cross-square section subjected to several explosive quantities (a-d) at $10$ m (\texttt{Method C}).}
	\label{f:meth_C}
\end{figure*}

We present in Figure \ref{f:meth_C} the prototype and (upscaled) model responses. For an explosive charge of $250$ kg of TNT equivalent (Fig. \ref{f:meth_C1}), the model predictions strongly differ from those of the prototype. While, for larger explosive weights (Fig.s \ref{f:meth_C2}-\ref{f:meth_C4}), the predictions agree remarkably well, especially for first peak response. It is worth noticing that the scaling laws do not account neither for wobbling motion nor for asynchronous arrival times of the shock wave at different points of the structure. As a result the equivalent impulse acting on the prototype structure is, in principle, not similar with that of the model. Nevertheless, for large (enough) charge weights, the similarity is preserved. Indeed, in these cases, rocking motion is dominant and the shock wave velocity is fast enough to have any important effect in the dynamic response of both the prototype and the model.

\bibliography{Biblio_all_clean.bib}

\begin{thebibliography}{}

\bibitem[ABAQUS, 2018]{abaqus}
ABAQUS (2018).
\newblock Abaqus analysis user's guide.
\newblock Technical Report Abaqus 6.14 Documentation, Simulia Corp.

\bibitem[Abou-Zeid et~al., 2011]{abou}
Abou-Zeid, B.~M., El-Dakhakhni, W.~W., Razaqpur, A.~G., and Foo, S. (2011).
\newblock {Response of Arching Unreinforced Concrete Masonry Walls to Blast
  Loading}.
\newblock {\em Journal of Structural Engineering}, 137(10):1205--1214.

\bibitem[Çaktı et~al., 2016]{CAKTI2016224}
Çaktı, E., Özden Saygılı, Lemos, J.~V., and Oliveira, C.~S. (2016).
\newblock Discrete element modeling of a scaled masonry structure and its
  validation.
\newblock {\em Engineering Structures}, 126:224 -- 236.

\bibitem[Baker et~al., 1991]{199173}
Baker, W.~E., Westine, P.~S., and Dodge, F.~T. (1991).
\newblock Simulating rigid body motion.
\newblock In {\em Similarity Methods in Engineering Dynamics}, volume~12 of
  {\em Fundamental Studies in Engineering}, pages 73 -- 95. Elsevier.

\bibitem[Bertrand, 1878]{bertrand1878homogeneite}
Bertrand, J. (1878).
\newblock Sur l'homog{\'e}n{\'e}it{\'e} dans les formules de physique.
\newblock {\em Cahiers de recherche de l'Academie de Sciences}, 86:916--920.

\bibitem[Casapulla et~al., 2017]{lourenco1}
Casapulla, C., Giresini, L., and Louren{\c{c}}o, P.~B. (2017).
\newblock Rocking and kinematic approaches for rigid block analysis of masonry
  walls: state of the art and recent developments.
\newblock {\em Buildings}, 7(3):69.

\bibitem[Cranz, 1925]{cranz}
Cranz, C. (1925).
\newblock Lehrbuch der ballistik.
\newblock {\em Julius Springer, Berlin}, 27.

\bibitem[Dimitrakopoulos and DeJong, 2012]{dejong}
Dimitrakopoulos, E.~G. and DeJong, M.~J. (2012).
\newblock Revisiting the rocking block: closed-form solutions and similarity
  laws.
\newblock {\em Proceedings of the Royal Society of London A: Mathematical,
  Physical and Engineering Sciences}, 468(2144):2294--2318.

\bibitem[Dragani{\'c} et~al., 2018]{draganic2018overview}
Dragani{\'c}, H., Varevac, D., and Luki{\'c}, S. (2018).
\newblock An overview of methods for blast load testing and devices for
  pressure measurement.
\newblock {\em Advances in Civil Engineering}, 2018.

\bibitem[Drosos and Anastasopoulos, 2014]{drosos2014shaking}
Drosos, V. and Anastasopoulos, I. (2014).
\newblock Shaking table testing of multidrum columns and portals.
\newblock {\em Earthquake Engineering \& Structural Dynamics},
  43(11):1703--1723.

\bibitem[Forg{\'a}cs et~al., 2017]{forgacs2017minimum}
Forg{\'a}cs, T., Sarhosis, V., and Bagi, K. (2017).
\newblock Minimum thickness of semi-circular skewed masonry arches.
\newblock {\em Engineering Structures}, 140:317--336.

\bibitem[Forg{\'a}cs et~al., 2018]{forgacs2018influence}
Forg{\'a}cs, T., Sarhosis, V., and Bagi, K. (2018).
\newblock Influence of construction method on the load bearing capacity of skew
  masonry arches.
\newblock {\em Engineering Structures}, 168:612--627.

\bibitem[Fragiadakis et~al., 2016]{fragiadakis2016vulnerability}
Fragiadakis, M., Stefanou, I., and Psycharis, I.~N. (2016).
\newblock Vulnerability assessment of damaged classical multidrum columns.
\newblock In {\em Computational Modeling of Masonry Structures Using the
  Discrete Element Method}, pages 235--253. IGI Global.

\bibitem[Friedlander, 1946]{friedlander}
Friedlander, F.~G. (1946).
\newblock The diffraction of sound pulses. i. diffraction by a semi-infinite
  plate.
\newblock {\em Proceedings of the Royal Society of London A}, 186:322 -- 344.

\bibitem[Gabrielsen et~al., 1975]{gabrielsen1975response}
Gabrielsen, B., Wilton, C., and Kaplan, K. (1975).
\newblock {Response of arching walls and debris from interior walls caused by
  blast loading}.
\newblock Technical report, URS Reasearch Company, San Mateo. CA.

\bibitem[Gilbert et~al., 2002]{gilbert2002performance}
Gilbert, M., Hobbs, B., and Molyneaux, T. (2002).
\newblock The performance of unreinforced masonry walls subjected to
  low-velocity impacts: experiments.
\newblock {\em International Journal of Impact Engineering}, 27(3):231--251.

\bibitem[Godio et~al., 2018]{Godio2018}
Godio, M., Stefanou, I., and Sab, K. (2018).
\newblock {Effects of the dilatancy of joints and of the size of the building
  blocks on the mechanical behavior of masonry structures}.
\newblock {\em Meccanica}, 53(7):1629--1643.

\bibitem[Hopkinson, 1915]{hopki}
Hopkinson, B. (1915).
\newblock British ordinance board minutes 13565.
\newblock {\em The National Archives, Kew, UK}, 11.

\bibitem[Kassotakis et~al., 2020]{kassotakis2020three}
Kassotakis, N., Sarhosis, V., Riveiro, B., Conde, B., D'Altri, A.~M., Mills,
  J., Milani, G., de~Miranda, S., and Castellazzi, G. (2020).
\newblock Three-dimensional discrete element modelling of rubble masonry
  structures from dense point clouds.
\newblock {\em Automation in Construction}, 119:103365.

\bibitem[Keys and Clubley, 2017]{KEYS2017229}
Keys, R.~A. and Clubley, S.~K. (2017).
\newblock {Experimental analysis of debris distribution of masonry panels
  subjected to long duration blast loading}.
\newblock {\em Engineering Structures}, 130:229 -- 241.

\bibitem[Kingery and Bulmash, 1984]{kingery}
Kingery, C.~N. and Bulmash, G. (1984).
\newblock {Technical report ARBRL-TR-02555: Air blast parameters from TNT
  spherical air burst and hemispherical burst}.
\newblock Technical report, U.S. Army Ballistic Research Laboratory.

\bibitem[Konstantinidis and N., 2005]{drumcolumns}
Konstantinidis, D. and N., M. (2005).
\newblock Seismic response analysis of multidrum classical columns.
\newblock {\em Earthquake Engineering \& Structural Dynamics},
  34(10):1243--1270.

\bibitem[Krauthammer and Altenberg, 2000]{KRAUTHAMMER20001}
Krauthammer, T. and Altenberg, A. (2000).
\newblock Negative phase blast effects on glass panels.
\newblock {\em International Journal of Impact Engineering}, 24(1):1 -- 17.

\bibitem[Li et~al., 2017]{LI2017107}
Li, Z., Chen, L., Fang, Q., Hao, H., Zhang, Y., Xiang, H., Chen, W., Yang, S.,
  and Bao, Q. (2017).
\newblock {Experimental and numerical study of unreinforced clay brick masonry
  walls subjected to vented gas explosions}.
\newblock {\em International Journal of Impact Engineering}, 104:107 -- 126.

\bibitem[Makris and Vassiliou, 2013]{makris2013planar}
Makris, N. and Vassiliou, M.~F. (2013).
\newblock Planar rocking response and stability analysis of an array of
  free-standing columns capped with a freely supported rigid beam.
\newblock {\em Earthquake Engineering \& Structural Dynamics}, 42(3):431--449.

\bibitem[Masi, 2020]{masithesis}
Masi, F. (2020).
\newblock {\em Fast-dynamic response and failure modes of masonry structures of
  non-standard geometry subjected to blast loads}.
\newblock PhD thesis, École Centrale de Nantes, Nantes, France.

\bibitem[Masi et~al., 2020a]{masi2020discrete}
Masi, F., Stefanou, I., Maffi-Berthier, V., and Vannucci, P. (2020a).
\newblock A discrete element method based-approach for arched masonry
  structures under blast loads.
\newblock {\em Engineering Structures}, 216:110721.

\bibitem[Masi et~al., 2019a]{hstamwall}
Masi, F., Stefanou, I., Vannucci, P., and Maffi-Berthier, V. (2019a).
\newblock {A Discrete Element Method approach for the preservation of the
  architectural heritage against explosions}.
\newblock In {\em Proc. of the 12th HSTAM International Congress on Mechanics},
  Thessaloniki, Greece.

\bibitem[Masi et~al., 2019b]{rockingmasi}
Masi, F., Stefanou, I., Vannucci, P., and Maffi-Berthier, V. (2019b).
\newblock Rocking response of inverted pendulum structures under blast loading.
\newblock {\em International Journal of Mechanical Sciences}, 157-158:833 --
  848.

\bibitem[Masi et~al., 2020b]{masi2020resistance}
Masi, F., Stefanou, I., Vannucci, P., and Maffi-Berthier, V. (2020b).
\newblock Resistance of museum artefacts against blast loading.
\newblock {\em Journal of Cultural Heritage}.

\bibitem[Michaloudis and Gebbeken, 2019]{gebbeken}
Michaloudis, G. and Gebbeken, N. (2019).
\newblock {Modeling masonry walls under far-field and contact detonations}.
\newblock {\em International Journal of Impact Engineering}, 123:84 -- 97.

\bibitem[Neils, 2005]{neils2005parthenon}
Neils, J. (2005).
\newblock {\em The Parthenon: from antiquity to the present}.
\newblock Cambridge University Press.

\bibitem[Pe\~{n}a et~al., 2007]{pena}
Pe\~{n}a, F., Prieto, F., Lourenço, P.~B., Campos~Costa, A., and Lemos, J.~V.
  (2007).
\newblock On the dynamics of rocking motion of single rigid$‐$block
  structures.
\newblock {\em Earthquake Engineering \& Structural Dynamics},
  36(15):2383--2399.

\bibitem[Psycharis et~al., 2000]{psycharis}
Psycharis, I., Papastamatiou, D.~Y., and Alexandris, A.~P. (2000).
\newblock Parametric investigation of the stability of classical columns under
  harmonic and earthquake excitations.
\newblock {\em {Earthquake Engineering and Structural Dynamics}},
  29:1093--1109.

\bibitem[Sarhosis et~al., 2016]{sarhosis2016stability}
Sarhosis, V., Asteris, P., Wang, T., Hu, W., and Han, Y. (2016).
\newblock On the stability of colonnade structural systems under static and
  dynamic loading conditions.
\newblock {\em Bulletin of Earthquake Engineering}, 14(4):1131--1152.

\bibitem[Stefanou et~al., 2015]{stefanou2015seismic}
Stefanou, I., Fragiadakis, M., and Psycharis, I.~N. (2015).
\newblock Seismic reliability assessment of classical columns subjected to near
  source ground motions.
\newblock In {\em Seismic assessment, behavior and retrofit of heritage
  buildings and monuments}, pages 61--82. Springer.

\bibitem[Stefanou et~al., 2011a]{STEFANOU20114325}
Stefanou, I., Psycharis, I., and Georgopoulos, I.-O. (2011a).
\newblock {Dynamic response of reinforced masonry columns in classical
  monuments}.
\newblock {\em Construction and Building Materials}, 25(12):4325 -- 4337.

\bibitem[Stefanou et~al., 2011b]{stefanou2011dynamic}
Stefanou, I., Vardoulakis, I., and Mavraganis, A. (2011b).
\newblock Dynamic motion of a conical frustum over a rough horizontal plane.
\newblock {\em International Journal of Non-Linear Mechanics}, 46(1):114--124.

\bibitem[Stockdale et~al., 2020]{stockdale2020seismic}
Stockdale, G.~L., Sarhosis, V., and Milani, G. (2020).
\newblock Seismic capacity and multi-mechanism analysis for dry-stack masonry
  arches subjected to hinge control.
\newblock {\em Bulletin of Earthquake Engineering}, 18(2):673--724.

\bibitem[USACE, 2008]{ufc08}
USACE (2008).
\newblock {UFC 3-340-02: Structures to Resist the Effects of Accidental
  Explosions}.
\newblock Technical report, U.S. Army.

\bibitem[Vannucci et~al., 2017a]{vannucci2017comparative}
Vannucci, P., Masi, F., and Stefanou, I. (2017a).
\newblock A comparative study on the effects of blast actions on a monumental
  structure.
\newblock Technical report, UVSQ ENPC.

\bibitem[Vannucci et~al., 2019]{VANNUCCI20192}
Vannucci, P., Masi, F., Stefanou, I., and Maffi-Berthier, V. (2019).
\newblock {Structural integrity of Notre Dame Cathedral after the fire of April
  15th, 2019}.
\newblock Technical report, UVSQ, ENPC, Ingérop.

\bibitem[Vannucci et~al., 2017b]{cnrs}
Vannucci, P., Stefanou, I., and Masi, F. (2017b).
\newblock {Report of the project "Cath{\'e}drales durables"}.
\newblock Technical report, CNRS, Paris, France.

\bibitem[Varma et~al., 1997]{varma}
Varma, R.~K., Tomar, C. P.~S., Parkash, S., and Sethi, V.~S. (1997).
\newblock {Damage to brick masonry panel walls under high explosive
  detonations}.
\newblock {\em Pressure vessels and piping division. ASME}, 351:207 -- 216.

\bibitem[Vassiliou, 2018]{vassiliou2018seismic}
Vassiliou, M.~F. (2018).
\newblock Seismic response of a wobbling 3d frame.
\newblock {\em Earthquake Engineering \& Structural Dynamics},
  47(5):1212--1228.

\bibitem[Voyagaki et~al., 2013]{elia}
Voyagaki, E., Psycharis, I.~N., and Mylonakis, G. (2013).
\newblock Rocking response and overturning criteria for free standing rigid
  blocks to single—lobe pulses.
\newblock {\em Soil Dynamics and Earthquake Engineering}, 46:85--95.

\bibitem[Wang et~al., 2012]{wang2012experimental}
Wang, W., Zhang, D., Lu, F., Wang, S.-C., and Tang, F. (2012).
\newblock Experimental study on scaling the explosion resistance of a one-way
  square reinforced concrete slab under a close-in blast loading.
\newblock {\em International Journal of Impact Engineering}, 49:158--164.

\bibitem[Zhang and Makris, 2001]{makris}
Zhang, J. and Makris, N. (2001).
\newblock Rocking response of free-standing blocks under cycloidal pulses.
\newblock {\em Journal of Engineering Mechanics}, 127(5):473--483.

\end{thebibliography}

\end{document}